
\input amstex
\documentstyle{amsppt}

\magnification=\magstep1
\hsize=6.5truein
\vsize=9truein

\font \smallrm=cmr10 at 10truept
\font \smallsl=cmsl10 at 10truept
\font \ssmallrm=cmr7 at 7truept
\font \smallbf=cmbx10 at 10truept
\font \smallit=cmti10 at 10truept

\baselineskip=.15truein

\def \longhookrightarrow {\lhook\joinrel\relbar\joinrel\rightarrow}
\def \llonghookrightarrow
{\lhook\joinrel\relbar\joinrel\relbar\joinrel\relbar\joinrel\rightarrow}

\def \smallcirc {\, {\scriptstyle \circ} \,}
\def \aij {a_{i j}}
\def \unon {1,\dots,n}
\def \N {\Bbb N}
\def \Z {\Bbb Z}
\def \C {\Bbb C}
\def \qm {q^{-1}}
\def \Ch {{\Bbb C}[[h]]}
\def \Cq {{\Bbb C}(q)}
\def \Cqqm {{\Bbb C} \! \left[ q, \qm \right]}
\def \rcalzero {{{\Cal R}^{(0)}}}
\def \rcaluno {{{\Cal R}^{(1)}}}

\def \runo {{R^{(1)}}}
\def \runoalpha {{R^{(1)}_{\, \alpha}}}
\def \rcalunoalpha {{{\Cal R}^{(1)}_\alpha}}
\def \gerb {\frak b}
\def \gerg {\frak g}

\def \uhg {U_h(\gerg)}
\def \uqQg {U_q^Q\!(\gerg)}
\def \uqMg {U_q^M\!(\gerg)}
\def \uqPg {U_q^P\!(\gerg)}
\def \caluqMg {{\Cal U}_q^M\!(\gerg)}
\def \caluqPg {{\Cal U}_q^P\!(\gerg)}

\def \caluepsilonPg {{\Cal U}_\varepsilon^P\!(\gerg)}
\def \caluunoMg {{\Cal U}_1^M\!(\gerg)}
\def \caluunoPg {{\Cal U}_1^P\!(\gerg)}

\def \uqQbm {U_q^Q({\frak b}_-)}

\def \uqQbp {U_q^Q({\frak b}_+)}
\def \ebar {\overline E}
\def \fbar {\overline F}


\document

\topmatter

{\ }

\vskip-49pt

   \hfill   {{\smallsl Communications in Mathematical Physics},  {\smallbf 184},
{\smallrm no.~1 (1997), 95--117}}   \hskip19pt   {\ }
                                            \par   
   \hfill   {{\smallbf DOI:}  {\smallrm 10.1007/s002200050054}}   \hskip31pt   {\ }
 \vskip4pt
   \hfill   {\smallrm {\smallsl The original publication is available at\/} \  www.springerlink.com/content/9c34gj7gm38xglta/}   \hskip19pt   {\ }

\vskip49pt

\title
    Geometrical meaning of  {\it R\/}-matrix  action for quantum groups at
roots of 1
\endtitle

\author
   Fabio Gavarini
\endauthor

\rightheadtext{ Geometrical meaning of  {\it R\/}-matrix  action for
quantum groups... }

\affil
   Dipartimento di Matematica, Istituto G. Castelnuovo,\\ Universit\`a
degli studi di Roma "La Sapienza" \\
\endaffil

\address\hskip-\parindent
   Dipartimento di Matematica, Istituto  "G. Castelnuovo"  \newline
   Universit\`a di Roma "La Sapienza"  \newline
   Piazzale Aldo Moro 5  \newline
   00185 Roma --- ITALY  \newline
   e-mail address:  \  gavarini\@mat.uniroma1.it
\endaddress

\abstract
    The present work splits in two parts: first, we perform a
straightforward generalization of results from [Re], proving
that quantum groups  $ \uqMg $  and their unrestricted
specializations at roots of 1,  in particular the function algebra  $ F[H]
$  of the Poisson group  $ H $  dual of  $ G $,  are braided;  second, as a main
contribution, we prove the convergence of the  (specialized)  $ R $-matrix
action to a birational automorphism of a  $ 2 \ell $-fold  ramified
covering of
$ {Spec \left( U_\varepsilon^M \! (\gerg) \right)}^{\times 2} $  when
$ \varepsilon $  is a primitive  $ \ell $-th  root of  $ 1 $,  and of a
$ 2 $-fold  ramified covering of  $ H $,  thus giving a geometric content to
the notion of triangularity (or braiding) for quantum groups at roots of 1.
 \endabstract

\endtopmatter

\vskip20pt

\centerline {\bf \S\; 1 \ Definitions }

\vskip10pt

   {\bf 1.1  Cartan data.}  \;  Let  $ \gerg $  be a complex finite
dimensional semisimple Lie algebra  of rank  $ n $,  with Cartan matrix  $
A:= {(\aij)}_{i,j=\unon} $;  let
$ R $  be its root system,  $ Q $,  resp.~$ P $,  be its root lattice,
resp.~weight lattice; we fix a  subset  $ R^+ \, (\subset R ) $  of
positive roots, a basis  $ \{ \alpha_1, \dots, \alpha_n \} \, (\subset R^+)
$  of simple roots, and we let  $ \{ \omega_1, \dots, \omega_n \} $  be the
dual basis of
$ P $.  We denote by  $ W $  the Weyl group of  $ \gerg $,  with
generators  $ s_1, \ldots, s_n $  (namely the reflections associated with
simple roots),  and we set $ \, N:= \#(R^+) \, $.  Finally, we let  $ (d_1,
\dots, d_n) $  be the (unique)  $ n $-tuple  of relatively prime positive
integers such that  $ \, (d_i \aij)_{i,j=\unon} \, $  is a symmetric
positive definite matrix.

\vskip7pt

{\bf 1.1  Quantum enveloping algebras.} \;  We briefly recall some
definitions.  The quantized universal enveloping algebra  $ \uhg $  is the
associative algebra with 1 over
$ \Ch $  generated by  $ \, Y_1, \dots , Y_n, H_1, \dots , H_n, X_1, \dots,
X_n \, $  with relations (for  $ i, j= \unon $)
  $$  \displaylines {
   H_i H_j - H_j H_i = 0  \cr
   H_i X_j - X_j H_i = a_{ij} X_j   \, , \qquad H_i Y_j - Y_j H_i = -a_{ij}
Y_j  \cr
   X_i Y_j - Y_j X_i = \delta_{ij} {\exp(h d_i H_i) - \exp(-h d_i H_i)
\over \exp(h d_i) - \exp(-h d_i)}  \cr }  $$
  $$  \eqalignno{
    \sum\limits_{k = 0}^{1-a_{ij}} (-1)^k {\left[ { 1-a_{ij} \atop k }
\right]}_{q_i} X_i^{1-a_{ij}-k} X_j X_i^k  &  = 0  &   \forall \; i \ne j
\cr
    \sum\limits_{k = 0}^{1-a_{ij}} (-1)^k {\left[ { 1-a_{ij} \atop k }
\right]}_{q_i} Y_i^{1-a_{ij}-k} Y_j Y_i^k  &  = 0 &   \forall \; i \ne j
\cr  }  $$ where  $ \, q:= \exp(h) \, $,  $ \, q_i:= q^{d_i}=\exp(h d_i) \,
$,  and the Gaussian binomial  $ \, {\left[{m \atop n} \right]}_q \, $  is
defined by
  $$  {\left[{m \atop n} \right]}_q := {{[m]}_q ! \over {[m-n]}_q ! {[n]}_q
!} \, ,  \qquad {[k]}_q ! := \prod_{s=1}^k {[s]}_q \, ,  \qquad {[s]}_q :=
{q^s - q^{-s} \over q - q^{-1}}  $$ for all  $ m $,  $ n $,  $ k $,  $ s
\in \N_+ $,  $ n \leq m $,  with  $ {[s]}_q , {[k]}_q ! , {\left[{m \atop
n} \right]}_q \in \Cqqm $.  It is known that  $ \, \uhg \, $  has a Hopf
algebra structure, given by ($ i= \unon $)
  $$  \matrix
   \Delta(Y_i) := Y_i \otimes \exp(-h d_i H_i) + 1 \otimes Y_i \, ,  &
S(Y_i):= - Y_i \exp(h d_i H_i) \, ,  &  \epsilon(Y_i) := 0 \, \phantom{.}
\\
   \Delta(H_i) := H_i \otimes 1 + 1 \otimes H_i \quad ,  &  S(H_i):= -H_i
\, ,  &
\epsilon(H_i) := 0 \, \phantom{.}  \\
 \Delta(X_i) := X_i \otimes 1 + \exp(h d_i H_i) \otimes X_i \, ,  &
S(X_i):= - \exp(-h d_i H_i) X_i \, ,  &  \epsilon(X_i) := 0 \, .  \\
      \endmatrix  $$
   \indent   Let  $ M $  be a lattice such that  $ \, Q \leq M \leq P \,
$:  the quantized  universal enveloping algebra  $ \uqMg $  (cf.~[DP], \S
9) is the associative algebra with 1 over  $ \Cq $  generated by  $ \, F_1,
\dots , F_n $,  $ L_\mu \, (\forall \, \mu \in M) $,
$ \, E_1, \dots, E_n \, $  with relations  ($ i, j= \unon; \, \mu, \nu \in
M $)
  $$  \displaylines {
   L_\mu L_\nu = L_{\mu + \nu} = L_\nu L_\mu  \quad , \qquad L_0 = 1  \cr
   L_{\mu} E_j = q^{\langle \mu \vert \alpha_j \rangle} E_j L_{\mu} \, ,
\qquad  L_{\mu} F_j = q^{- \langle \mu \vert \alpha_j \rangle} F_j L_{\mu}
\cr
   E_i F_j - F_j E_i = \delta_{ij} {L_{\alpha_i} - L_{-\alpha_i} \over q_i
- q_i^{-1}}  \cr }  $$
  $$  \eqalignno{
    \sum_{k = 0}^{1-a_{ij}} (-1)^k {\left[ { 1-a_{ij} \atop k }
\right]}_{q_i} E_i^{1-a_{ij}-k} E_j E_i^k  &  = 0  &   \forall \; i \ne j
\cr
    \sum_{k = 0}^{1-a_{ij}} (-1)^k {\left[{1-a_{ij} \atop k}
\right]}_{q_i} F_i^{1-a_{ij}-k} F_j F_i^k  &  = 0  &   \forall \; i \ne j
\cr }  $$
   \indent   A Hopf algebra structure on  $ \, \uqMg \, $  is defined by
($ i= \unon; \, \mu \in M $)
  $$  \matrix
   \Delta(F_i) := F_i \otimes L_{-\alpha_i} + 1 \otimes F_i \,,  &
S(F_i):= -F_i L_{\alpha_i} \, ,  &  \epsilon(F_i) := 0 \, \phantom{.}  \\
   \Delta(L_\mu) := L_\mu \otimes L_\mu \, ,  &   S(L_\mu):= L_{-\mu} \, ,
&  \epsilon(L_\mu) := 1 \, \phantom{.}  \\
  \Delta(E_i) := E_i \otimes 1 + L_{\alpha_i} \otimes E_i \, ,  &
S(E_i):= -L_{-\alpha_i} E_i \, ,  &  \epsilon(E_i) := 0 \, .  \\
      \endmatrix  $$
   \indent   It is clear that  $ \, \uqMg \hookrightarrow U_q^{M'} \!
(\gerg) \, $  whenever  $ \, Q \leq M \leq M' \leq P \, $,  this being a
Hopf algebra embedding.  In the sequel we shall also use notation  $ \, L_i
:= L_{\omega_i} $,  $ K_i := L_{\alpha_i} \, $  (for all  $ i= \unon $).
                                                  \par
   The very definitions imply the existence of a Hopf algebra monomorphism
  $$  j \colon \uqQg \llonghookrightarrow \uhg  $$ given by  $ \; q \mapsto
\exp(h) \, $,  $ \, F_i \mapsto Y_i \, $,  $ \, K_i \mapsto \exp(h d_i H_i)
\, $,  $ \, E_i \mapsto X_i \, $;  still from definitions it is also clear
that this uniquely extends to an embedding
  $$  j \colon \uqMg \llonghookrightarrow \uhg  $$ for all lattices  $ M
$;  in particular  $ \, \uqPg \longhookrightarrow \uhg \, $.  Finally, we
shall denote by  $ U_q^M\!(\gerb_+) $,  resp.~$ U_q^M\!(\gerb_-) $,  the
Hopf subalgebra of  $ \uqMg $  generated by  $ L_\mu $'s  and
$ E_i $'s,  resp.~$ F_i $'s:  these are called  {\it quantum Borel
subalgebras}.

\vskip 7pt

   An interesting property that Hopf algebras can enjoy is
quasitriangularity:

\vskip7pt

\proclaim {Definition 1.2}  (cf.~[Dr]) A Hopf algebra  $ H $  is called
quasitriangular if there exists an invertible element  $ R \in H \otimes H
$  (or an element of an appropriate completion of  $ H \otimes H $)  such
that
  $$  R \cdot \Delta (a) \cdot R^{-1} = {\hbox{\rm Ad}}(R) (\Delta (a)) =
\Delta^{\hbox{\smallrm op}}(a)
\eqno (1.1)  $$
  $$  (\Delta \otimes {\hbox{\rm id}}) (R) = R_{13} R_{23}    \eqno (1.2)
$$
  $$  ({\hbox{\rm id}} \otimes \Delta) (R) = R_{13} R_{12}    \eqno (1.3)
$$ where  $ \, \Delta^{\hbox{\smallrm op}} \, $  is the opposite
comultiplication, i.~e.  $ \,
\Delta^{\hbox{\smallrm op}}(a) = \sigma \smallcirc \Delta(a) \, $  with  $
\, \sigma \colon
\, A^{\otimes 2} \to A^{\otimes 2} \, $,  $ \, a \otimes b \mapsto b
\otimes a  \, $,  and  $ \, R_{12}, R_{13}, R_{23} \in H^{\otimes 3} \, $
(or the appropriate completion of  $ H^{\otimes 3} $),  $ \, R_{12} = R
\otimes 1 \, $,  $ \, R_{23} = 1 \otimes R \, $,  $ \, R_{13} = (\sigma
\otimes {\hbox{\rm id}}) (R_{23}) =  ({\hbox{\rm id}} \otimes \sigma)
(R_{12}) \, $.   $ \square $
\endproclaim

\vskip7pt

   As a corollary of this definition,  $ R $  satisfies the Yang-Baxter
equation in  $ H^{\otimes 3} $  (cf.~[Ta])
  $$  R_{12} R_{13} R_{23} = R_{23} R_{13} R_{12}  $$
so that a braid group action is defined on tensor products of  $ H $--modules.
The quantum universal enveloping algebra  $ \uhg $  is quasitriangular
(cf.~[Dr], [LS], [KR]):  its  $ R $-matrix  is
  $$  R = \prod_{\alpha \in R^+} \exp_{q_\alpha} \! \Big( \left(
{q_\alpha}^{-1} - q_\alpha \right)  X_\alpha \otimes Y_\alpha \Big) \cdot
\exp \! \bigg( -h \sum_{i,j=1}^n B_{ij} H_i \otimes H_j \bigg)  $$  where
$ \prod_{\alpha \in R^+} $  denotes an ordered product (with respect to a
fixed convex ordering of  $ R^+ \, $),  $ \, q_\alpha :=  q^{d_\alpha}
\, $  (where  $ d_\alpha $  is one-half the square length of the root
$ \alpha $; in particular  $ \, d_{\alpha_i} = d_i \, $  for all  $ i $),
$ \, {(B_{ij})}_{i,j = \unon} :={(d_i \aij)}_{i,j=\unon}^{-1} \, $,  and  $
X_\alpha $,
$ Y_\alpha $  are  $ q $-analogue of root vectors (not unique, however)
attached to roots  $ \alpha $,
$ -\alpha $.  On the other hand, this is not true   --- strictly speaking
---  for the  $ \Cq $--algebras
$ \uqMg $:  to be precise we need a slight modification of the notion of
quasitriangularity, suggested by Reshetikin, as follows:

\vskip7pt

\proclaim {Definition 1.3}  (cf. [Re], Definition 2) A Hopf algebra  $ H $
is called braided if there exists an automorphism  $ {\Cal R}
$  of  $ H \otimes H $  (or of an appropriate completion of  $ H \otimes H
$)  distinct from  $ \sigma \colon a \otimes b \mapsto b \otimes a $  such
that
  $$  {\Cal R} \smallcirc \Delta = \Delta^{\hbox{\smallrm op}}
\eqno  (1.4)  $$
  $$  (\Delta \otimes {\hbox{\rm id}}) \smallcirc {\Cal R} = {\Cal R}_{13}
\smallcirc {\Cal R}_{23}
\smallcirc (\Delta \otimes{\hbox{\rm id}})   \eqno  (1.5)  $$
  $$  ({\hbox{\rm id}} \otimes \Delta) \smallcirc {\Cal R} = {\Cal R}_{13}
\smallcirc {\Cal R}_{12}
\smallcirc ({\hbox{\rm id}} \otimes \Delta)    \eqno  (1.6)  $$
   \indent   Here  $ {\Cal R}_{12}, {\Cal R}_{13}, {\Cal R}_{23} $  are the
automorphisms of
$ H \otimes H \otimes H $  defined by  $ \, {\Cal R}_{12} = {\Cal R}
\otimes {\hbox{\rm id}} \, $,
$ \, {\Cal R}_{23} = {\hbox{\rm id}} \otimes {\Cal R} \, $,  $ \, {\Cal
R}_{13} = (\sigma \otimes {\hbox{\rm id}}) \smallcirc ({\hbox{\rm id}}
\otimes {\Cal R}) \smallcirc (\sigma \otimes {\hbox{\rm id}}) \, $.   $
\square $
\endproclaim

\vskip7pt

   It follows from this definition that  $ {\Cal R} $  satisfies the
Yang-Baxter equation in  $ End(H^{\otimes 3}) $:
  $$  {\Cal R}_{12} \smallcirc {\Cal R}_{13} \smallcirc {\Cal R}_{23} =
{\Cal R}_{23} \smallcirc {\Cal R}_{13} \smallcirc {\Cal R}_{12}
\eqno  (1.7)  $$
which yields a braid group action on tensor powers of  $ H $.  Furthermore,
it is clear that if  $ (H,R) $  is quasitriangular, then  $ \big( H,
{\hbox{\rm Ad}}(R) \big) $  is braided.  Again from [Re] we
resume another definition (slightly modified indeed).

\vskip7pt

\proclaim {Definition 1.4}  (cf. [Re], Definition 3) Let
$ H $  be a Hopf algebra, let  $ \rcalzero $  be an algebra automorphism
of  $ H \otimes H $  and  $ \runo \in H \otimes H $  an invertible element
such that
  $$  \runo \cdot \rcalzero \big( \Delta(a) \big) \cdot \runo^{-1} = \Big(
{\hbox{\rm Ad}} \big( \runo
\big) \smallcirc \rcalzero \Big) \big( \Delta(a) \big) =
\Delta^{\hbox{\smallrm op}}(a)   \eqno (1.8)  $$
  $$  (\Delta \otimes \hbox{\rm id}) \smallcirc \rcalzero = \rcalzero_{\!
13} \smallcirc \rcalzero_{\! 23} \smallcirc (\Delta \otimes \hbox{\rm
id})   \eqno  (1.9)  $$
  $$  (\hbox{\rm id} \otimes \Delta) \smallcirc \rcalzero = \rcalzero_{\!
13} \smallcirc \rcalzero_{\! 12} \smallcirc ({\hbox{\rm id}} \otimes
\Delta)    \eqno  (1.10)  $$
  $$  (\Delta \otimes \hbox{\rm id}) \big( \runo \big) = \runo_{\! 13}
\cdot {\rcalzero}_{\! 13} \big( \runo_{\! 23} \big)   \eqno (1.11)  $$
  $$  (\hbox{\rm id} \otimes \Delta) \big( \runo \big) = \runo_{\! 13}
\cdot {\rcalzero}_{\! 13} \big( \runo_{\! 12} \big)   \eqno (1.12)  $$
then  $ \Big( H, {\hbox{\rm Ad}} \big( \runo \big) \smallcirc \rcalzero
\Big) $  is a braided Hopf algebra, and the element  $ \runo $  is called the
{\it universal  $ R $-matrix}  of  $ \Big( H, {\hbox{\rm
Ad}} \big( \runo \big) \smallcirc \rcalzero \Big) \, $.
$ \square $
\endproclaim

\vskip7pt

   Finally we recall from [Ta] the strictly related notion below:

\vskip7pt

\proclaim {Definition 1.5}  (cf. [Ta], \S~4)  Let  $ H $  be a Hopf
algebra, let   $ \Phi $  be an algebra automorphism  of  $ H \otimes H $
and  $ C \in H \otimes H $  be an invertible element such that
  $$  C^{-1} \cdot \Phi(\Delta^{\hbox{\smallrm op}}(a)) \cdot C =
\Delta(a)   \eqno  (1.13)  $$
  $$  (\Phi_{23} \smallcirc \Phi_{13}) (C_{12}) = C_{12}         \eqno
(1.14)  $$
  $$  (\Phi_{12} \smallcirc \Phi_{13}) (C_{23}) = C_{23}         \eqno
(1.15)  $$
  $$  (\Delta \otimes {\hbox{\rm id}}) (C) = \Phi_{23} (C_{13}) \cdot
C_{23}   \eqno  (1.16)  $$
  $$  ({\hbox{\rm id}} \otimes \Delta) (C) = \Phi_{12} (C_{13}) \cdot
C_{12}   \eqno  (1.17)  $$ then we will say that  $ (H,C,\Phi) $  is a
pretriangular Hopf algebra.   $ \square $
\endproclaim

\eject

\centerline {\bf \S\; 2 \  Some  $ q $-calculus }

\vskip10pt

{\bf 2.1.} \;  In this section we introduce some material to be used in the
sequel; as  standard references for  $ q $-special  functions and related
matters we quote [Ex] and [GR].
                                                \par
   Let us introduce some  $ q $-symbols.  We have  $ q $-numbers
  $$  {(s)}_q := {\,q^s - 1\, \over \,q - 1\,} \, ,  \qquad {(k)}_q ! :=
\prod_{s=1}^k  {(s)}_q \, ,  \qquad {\left({m \atop n} \right)}_q :=
{\,{(m)}_q !\, \over \,{(m-n)}_q ! \,  {(n)}_q !\,}  $$ for  $ s, k, m, n
\in \N_+ $,  $ n \leq m \, $,  with
$ {(s)}_q , {(k)}_q ! , {\left({m \atop n} \right)}_q \in \Cqqm $,  and the
symbol
$ \, {(a;q)}_n := \prod_{k=0}^{n-1} (1 - a q^k) \, $,  for  $ n \in \N $,
$ a \in \C $.  Now consider the function of  $ z $
  $$  {(z;q)}_\infty := \prod_{n=0}^\infty
(1- z q^n)  \ ;  $$ we regard it as an element of  $ \Cq[[z]] $.  The
infinite  product expressing  $ {(z;q)}_\infty $  converges to an  analytic
function of  $ z $  in any finite part of  $ \C $   if  $ q $  is a complex
number such that  $ |q| < 1 $;  its  Taylor series is then
  $$  {(z;q)}_\infty = \sum_{n=0}^\infty {{(-1)}^n q^{n \choose 2} \over
{(q;q)}_n} z^n \; .  $$
   \indent   Also the series
  $$  e_q(z) := \sum_{n=0}^\infty {1 \over {(q;q)}_n} z^n \, ,  \quad
\qquad  E_q(z) := \sum_{n=0}^\infty {q^{n \choose 2} \over {(q;q)}_n} z^n
$$ both converge to analytic functions of  $ z $;  moreover, one has
  $$  e_q(z) = {{(z;q)}_\infty}^{\!\!\! -1} \; ,  \quad \qquad  E_q(z) =
{(-z;q)}_\infty  $$ so that  $ \, E_q(-z) e_q(z) = 1 \, $.  Finally
  $$  \exp_q(z) := \sum_{n=0}^\infty {\,1\, \over \, (n)_{q^2} ! \,} \,
z^n  $$ thus one has
  $$  \exp_q(z) = e_{q^2} \Big( \big( 1 - q^2 \big) z \Big) \, .  $$
   \indent   We claimed above that  $ {(z;q)}_\infty $  is an analytic
function of  $ z $  for  $ \, \vert q \vert < 1 \, $;  the following lemma
describes the behavior of this function for  $ \, q \rightarrow \varepsilon
\, $,   $ \, \varepsilon $  a root of 1.

\vskip7pt

\proclaim {Lemma 2.2}  ([Re], Lemma 3.4.1)  Let
$ \varepsilon $  be a primitive  $ \ell $-th  root of
$ 1 $,  with  $ \ell $  odd.  The asymptotic behavior of the function (of
$ q $)  $ {(z;q)}_\infty $  for  $ \, q \rightarrow \varepsilon \, $   is
given by
  $$  \eqalign{
   {(z;q)}_\infty  &  = \exp \left({-1 \over {q^{\ell^2} - 1}} \cdot
\int_0^{z^\ell}  {\log(1-t) \over t} \, dt \right) \cdot \prod_{k=0}^{\ell
- 1} {\left( 1 -
\varepsilon^k z \right)}^{k/\ell} \cdot \big( 1 + {\Cal O}(q - \varepsilon)
\big) =  \cr
   &  = \exp \left({1 \over {q^{\ell^2} - 1}} \sum_{n=1}^\infty {1 \over
n^2} \cdot z^{\ell n} \right) \cdot \prod_{k=0}^{\ell - 1} {\left( 1 -
\varepsilon^k z \right)}^{k/\ell} \cdot \big( 1 + {\Cal O}(q-\varepsilon)
\big) \; .  \cr }   \eqno (2.1)  $$
\endproclaim

\demo{Proof}  Taylor expansion of  $ \, \log(1-t) \, $  shows that the two
expressions in right-hand-side of (2.1) are equivalent.  Now, the function
$ \, {(z;q)}_\infty \, $  satisfies the difference equation
  $$  {\left( z q^\ell; q \right)}_\infty = {1 \over {(z;q)}_\ell} \,
{(z;q)}_\infty   \eqno (2.2)  $$
and it is uniquely determined by this
property along with the condition  $ \, (0;q)_\infty = 1 \, $.  But
  $$  \psi_z(q) :=  \exp\left({1 \over {q^{\ell^2} - 1}} \sum_{n=1}^\infty
{1 \over n^2} z^{\ell n} \right)
\cdot \prod_{k=0}^{\ell - 1} {(1 - \varepsilon^k z)}^{k/\ell}  $$
has the
asymptotic behavior, for  $ \, q \rightarrow \varepsilon \, $,  of  the
solution of (2.2); in fact we have
  $$  {\psi_z(q) \over {(z;q)}_\ell \cdot \psi_{z q^\ell}(q)} =  $$
  $$  = {\exp\left( {\left( {q^{\ell^2} - 1} \right)}^{-1} \cdot
\sum_{n=1}^\infty {1 \over n^2} z^{\ell n} \right) \cdot
\prod_{k=0}^{\ell - 1} {\left( 1 - \varepsilon^k z \right)}^{k/\ell}
\over  {\prod_{k=0}^{\ell - 1} \left( \! 1 - z q^k \right) \cdot \exp \! \left(
{\left( {q^{\ell^2} - 1} \right)}^{-1} \sum_{n=1}^\infty {1 \over n^2}
q^{\ell^{\scriptscriptstyle 2} n} z^{\ell n} \right) \cdot
\prod_{k=0}^{\ell - 1} {\left( \! 1 - \varepsilon^k z q^\ell
\right)}^{k/\ell} } } =  $$
  $$  = { \prod_{k=0}^{\ell - 1} {\left( 1 - \varepsilon^k z
\right)}^{k/\ell}  \over  {\prod_{k=0}^{\ell - 1} {\left( 1 - \varepsilon^k
z q^\ell \right)}^{k/\ell} \cdot
\prod_{k=0}^{\ell - 1} \left( 1 - z q^k \right)} } \cdot \exp \left(
{\left( {q^{\ell^2} - 1} \right)}^{-1} \cdot \sum_{n=1}^\infty {\left( 1 -
q^{\ell^2 n} \right) \over n^2} z^{\ell n} \right) =  $$
  $$  = \prod_{k=0}^{\ell - 1} {\left({1 - \varepsilon^k z  \over 1 -
\varepsilon^k z q^\ell}\right)}^{\! k/\ell} \cdot {1 \over
{\prod_{k=0}^{\ell - 1} \left( 1 - z q^k \right)}} \cdot
\exp\!\left( -\sum_{n=1}^\infty {{q^{\ell^2 n} - 1} \over {q^{\ell^2} - 1}}
\cdot {{\left( z^\ell \right)}^n \over n^2} \right) =  $$
  $$  = \prod_{k=0}^{\ell - 1} \! {\left({1 - \varepsilon^k z  \over  1 -
\varepsilon^k z q^\ell}\right)}^{\! k/\ell} \cdot {1 \over
{\prod_{k=0}^{\ell - 1} \left( 1 - z q^k \right)}} \cdot \exp\!\left(
-\sum_{n=1}^\infty (n)_{q^{\ell^{\scriptscriptstyle 2}}} {\left( z^\ell
\right)}^n \big/ n^2 \right) \; ;  $$
when  $ \, q \longrightarrow \varepsilon \, $  we have  $ \, \lim_{q \rightarrow
\varepsilon} (n)_{q^{\ell^{\scriptscriptstyle 2}}} = n \, $,
$ \lim_{q \rightarrow \varepsilon} \prod_{k=0}^{\ell - 1} {\left( {1 -
\varepsilon^k z  \over  1 - \varepsilon^k z q^\ell} \right)}^{\! k/\ell} =
1 \, $,                and\break
$ \, \lim_{q \rightarrow \varepsilon} {1  \over  \prod_{k=0}^{\ell
-1} \left( 1 - z q^k \right)} = {\left( \prod_{k=0}^{\ell -1} \left( 1 - z
\varepsilon^k \right) \right)}^{\! -1} = \prod_{k=0}^{\ell - 1}
\varepsilon^k \cdot \prod_{k=0}^{\ell -1} \left( \varepsilon^k - z \right) = 1 - z^\ell
\, $;  thus
  $$  \lim_{q \rightarrow \varepsilon} {\psi_z(q) \over {(z;q)}_\ell \cdot
\psi_{z q^\ell}(q)} = {1 \over 1 - z^\ell} \cdot \exp\left(
-\sum_{n=1}^\infty {1 \over n} {\left( z^\ell \right)}^n \right) = {\exp
\Big( \log \big(1 - z^\ell \big) \Big) \over 1 - z^\ell} = {1 - z^\ell
\over 1 - z^\ell} = 1 \; ,  $$ i.~e.~$ \, \lim_{q \rightarrow 1} {\psi_z(q)
\over {(z;q)}_\ell \cdot \psi_{z q^\ell}(q)} = 1 \, $.  Moreover from
definition  $ \, \psi_0(q) = 1 \, $.  The claim follows.   $ \square $
\enddemo
\eject

\vskip2,5truecm

\centerline {\bf \S\; 3 \  Braiding of quantum enveloping algebras }

\vskip10pt

{\bf 3.1.} \;  As we said, it is well known that quantum algebras  $ \uhg
$  are quasitriangular; this is proved by means of Drinfeld's method of the
"quantum double" (cf. [Dr] and others).  On the other hand, for the
$ \Cq $--algebras  $ \uqMg $  the correct statement is that they are braided;
for  $ \, \gerg = {\frak s}{\frak l}(2) \, $,  this is proved in [Re]: here
we quickly perform the (straightforward) generalization.
                                                 \par
  To begin with we define a suitable completion of
$ {\uqMg}^{\otimes 2} $,  namely
  $$  \uqMg^{\widehat \otimes 2} := \left\{\, \sum_{n=0}^{+\infty} {\Cal
E}_n \cdot P_n^- \otimes P_n^+ \cdot {\Cal F}_n \,\right\}  $$ where  $ \,
P_n^- \in U_q^M\!(\gerb_-) $,  $ P_n^+ \in U_q^M\!(\gerb_+) $  ($
U_q^M\!(\gerb_\pm) $  being opposite quantum Borel subalgebras),  $ {\Cal
E}_n \in \sum_{\vert \beta \vert = n} {\left( \uqMg \right)}_\beta $,  $
{\Cal F}_n \in \sum_{\vert \beta \vert = -n} {\left( \uqMg \right)}_\beta
\, $.  It is clear that  $ {\uqMg}^{\widehat \otimes 2} $  is a completion
of  $ {\uqMg}^{\otimes 2} $  as Hopf algebra.  From now on, as in [DD],
[DP], we set  $ E_\alpha := X_\alpha $,  $ F_\alpha := Y_\alpha $.

\vskip7pt

\proclaim {Theorem 3.2}  Let  $ \rcalzero $  be the algebra automorphism
of  $ {\uqMg}^{\widehat \otimes 2} $  defined by
  $$  \aligned
     \rcalzero (L_\mu \otimes 1) := L_\mu \otimes 1 \; ,  &  \quad
\rcalzero (1 \otimes L_\mu) := 1 \otimes L_\mu  \\
     \rcalzero(E_i \otimes 1) := E_i \otimes L_{-\alpha_i} \; ,  &  \quad
\rcalzero (1 \otimes E_i) := L_{-\alpha_i} \otimes E_i  \\
     \rcalzero (F_i \otimes 1) := F_i \otimes L_{\alpha_i} \; ,  &  \quad
\rcalzero (1 \otimes F_i) := L_{\alpha_i} \otimes F_i  \\
      \endaligned  $$ ($ i=\unon $;  $ \mu \in M $)  and let  $ \, \runo
\in {\uqMg}^{\widehat \otimes 2} \, $  be defined by
  $$  \runo := \prod_{\alpha \in R^+} \exp_{q_\alpha} \! \Big( \big(
{q_\alpha}^{-1} - q_\alpha \big) E_\alpha \otimes F_\alpha \Big)  $$
   \indent  Then  $ \, \left( \uqMg, {\hbox{\rm Ad}}(\runo) \smallcirc
\rcalzero \right) \, $  is a braided Hopf algebra, with  $ \runo $  as
$ R $-matrix.
\endproclaim

\demo{Proof}  We just outline the main steps, details being trivial.  First
of all, direct computation on generators shows that (1.9) and (1.10) hold.
Then define
$ \, C \in \uqQg^{\widehat \otimes 2} \left( \subset \uqMg^{\widehat
\otimes 2} \right) \, $  by
  $$  C := \sum_{\beta \in Q_+} q^{(\beta,\beta)} \cdot \left( K_\beta^{-1}
\otimes K_\beta \right) \cdot C_\beta  $$ where  $ \, Q_+ := Q \cap P_+ \,
$  is the positive root lattice and  $ C_\beta $  is the canonical element
of the bilinear pairing  $ \, {\left( \uqQbp \right)}_\beta \times {\left(
\uqQbm \right)}_{-\beta} \longrightarrow \Cq \, $  among quantum Borel
algebras; let also  $ \, \Phi:= \rcalzero^{-1} \, $;  then it is proved in
[Ta], Theorem 4.3.3 that  $ \, \left( \uqQg, C, \Phi \right) \, $  is a
pretriangular Hopf algebra; the same proof also works for
$ \uqMg $  instead of  $ \uqQg $.  Now trivial checking yields  $ \, \runo
= \Phi^{-1} (C) \, $  (using Proposition 3.7 in [DD]); therefore
$ \left( \uqMg, C, \Phi \right) $  being pretriangular implies that
$ \left( \uqMg, {\hbox{\rm Ad}}(\runo) \smallcirc \rcalzero \right) $  is
braided.   $ \square $
\enddemo

\vskip7pt

   {\bf Remark 3.3.} \;  Applying the remarks in \S 2 we can provide a
multiplicative formula for the universal  $ R $-matrix  $ \runo $  of  $
\uqMg $,  namely
  $$  \runo = \prod_{\alpha \in R^+} \exp_{q_\alpha}
\!\left( \left( {q_\alpha}^{-1} - q_\alpha \right) \cdot E_\alpha \otimes
F_\alpha \right) =   \hfill {\ }  $$
  $$  {\ } \hfill   = \prod_{\alpha \in R^+} e_{q_\alpha^2}
\!\left( \left( {q_\alpha}^{-1} - q_\alpha \right) \cdot \left( 1 -
q_\alpha^2 \right) \cdot E_\alpha \otimes F_\alpha \right) =  $$
  $$  = \prod_{\alpha \in R^+} e_{q_\alpha^2} \!\left( \left(
{q_\alpha}^{-1} - q_\alpha \right) \cdot q_\alpha \left( q_\alpha^{-1} -
q_\alpha \right) \cdot E_\alpha \otimes F_\alpha \right) =    \hfill {\ }
$$
  $$  {\ } \hfill    = \prod_{\alpha \in R^+} e_{q_\alpha^2} \left(
q_\alpha \ebar_\alpha \otimes
\fbar_\alpha \right) = \prod_{\alpha \in R^+} {{\left( q_\alpha \cdot
\ebar_\alpha \otimes \fbar_\alpha;
\, q_\alpha^2 \right)}_\infty}^{\!\!\! -1}  $$  where  $ \, \ebar_\alpha :=
\left( q_\alpha - q_\alpha^{-1} \right) E_\alpha \, $  and  $ \,
\fbar_\alpha := \left( q_\alpha - q_\alpha^{-1} \right) F_\alpha \, $
denote  {\sl modified}  root vectors; in other words
  $$  \runo = \prod_{\alpha \in R^+}  {{\left( q_\alpha \cdot \ebar_\alpha
\otimes \fbar_\alpha;
\, q_\alpha^2 \right)}_\infty}^{\!\!\! -1} \; .   \eqno (3.1)  $$

\vskip7pt

\proclaim{Definition 3.4} \;  We let  $ \caluqMg $  be the  $ \Cqqm
$--subalgebra  of  $ \uqMg $  generated by
  $$  \Big\{\, \fbar_\alpha, \, L_\mu, \, \ebar_\alpha \,\Big\vert\, \alpha
\in R^+, \mu \in M \,\Big\} \, .  $$
   \indent   Furthermore, for any  $ \, c \in \C \, $  we let
  $$  {\Cal U}_c^M \! (\gerg) := \, \caluqMg \Big/ (q-c) \,
\caluqMg \, \cong \, \caluqMg \otimes_{\C \left[ q, \qm \right]} \C  $$
(with  $ \, \C \cong \Cqqm \big/ (q-c) \, $)  be the corresponding
specialized algebra.   $ \square $
\endproclaim

\vskip7pt

   {\bf Remark 3.5.} \;  The previous definition is different but
equivalent to the original one in [DP], \S 12, equivalence arising from the
very description of
$ \caluqMg $  made therein.  It is also proved in [DP] that  $ \caluqMg $
is a  $ \Cqqm $--integer  form of  $ \uqMg $.

\vskip7pt

   {\bf 3.6.} \;  Our goal now is to show that  $ \caluqMg $  is
braided: to be precise, we  could say that the
braiding structure of  $ \uqMg $  gives by restriction a
braiding structure for  $ \caluqMg $.  To begin with, we define
a suitable completion of
$ \caluqMg^{\otimes 2} $  (mimicking \S 3.1),  namely
  $$  \caluqMg^{\widehat \otimes 2} := \left\{\,
\sum_{n=0}^{+\infty} {\overline{\Cal E}}_n \cdot P_n^-
\otimes P_n^+ \cdot {\overline{\Cal F}}_n \,\right\}  $$ where  $ \, P_n^-
\in {\Cal U}_q^M\!(\gerb_-) $,  $ P_n^+ \in {\Cal U}_q^M\!(\gerb_+) $,
$ {\overline{\Cal E}}_n \in \sum_{\vert \beta \vert = n} {\left( \caluqMg
\right)}_\beta $,
$ {\overline{\Cal F}}_n \in \sum_{\vert \beta \vert = -n} {\left( \caluqMg
\right)}_\beta \, $.  It is clear that  $ \caluqMg^{\widehat \otimes 2} $
is a completion of  $ \caluqMg^{\otimes 2} $  as Hopf algebra; moreover we
have  $ \, \caluqMg^{\widehat \otimes 2} \subseteq \uqMg^{\widehat \otimes
2} \, $  via the natural embedding  $ \, \caluqMg \longhookrightarrow \uqMg
\, $.

\vskip7pt

\proclaim {Proposition 3.7}  The restriction of
$ \rcalzero $  (cf.~Corollary 3.3) to  $ \caluqMg^{\widehat \otimes 2} $  is
defined by
  $$  \aligned
     \widetilde\rcalzero (L_\mu \otimes 1) := L_\mu \otimes 1 \; ,  &
\quad
\widetilde\rcalzero (1 \otimes L_\mu) := 1 \otimes L_\mu   \\
     \widetilde\rcalzero(\ebar_\alpha \otimes 1) := \ebar_\alpha \otimes
L_{-\alpha} \; ,  &  \quad   \widetilde\rcalzero (1 \otimes \ebar_\alpha) :=
L_{-\alpha} \otimes \ebar_i   \\
     \widetilde\rcalzero (\fbar_\alpha \otimes 1) := \fbar_\alpha \otimes
L_\alpha
\; ,  &  \quad   \widetilde\rcalzero (1 \otimes \fbar_\alpha) := L_\alpha
\otimes
\fbar_\alpha   \\
       \endaligned  $$ ($ \, \mu \in M, \, \alpha \in R^+ \, $)  so that  $
\rcalzero $  restricts to an algebra automorphism
$ \widetilde\rcalzero $  of
$ \caluqMg^{\widehat \otimes 2} $.  Moreover,  let  $ \, \runo \in
\uqMg^{\widehat \otimes 2} \, $  be defined (as in Corollary 3.3) by
  $$  \runo := \prod_{\alpha \in R^+} \exp_{q_\alpha} \!\Big( \left(
{q_\alpha}^{-1} - q_\alpha \right)  E_\alpha \otimes F_\alpha \Big) =
\prod_{\alpha \in R^+} {{\left( q_\alpha \cdot \ebar_\alpha \otimes
\fbar_\alpha; \, q_\alpha^2 \right)}_\infty}^{\!\!\! -1} \; .  $$
   \indent   Then  $ \hbox{\rm Ad} \left( \runo \right) $  restricts to an
automorphism
$ \widetilde\rcaluno $  of  $ \caluqMg^{\widehat \otimes 2} $,  and  $ \,
\left( \caluqMg, \widetilde{\Cal R} \right) \, $   --- with  $ \,
\widetilde{\Cal R} :=
\widetilde\rcaluno \smallcirc \widetilde\rcalzero \, $ ---   is a
braided Hopf algebra.
\endproclaim

\demo{Proof}  The first part of the statement is trivial.  As for the
second, we must recall that the specialization
$ \, \caluunoMg := \, \caluqMg \Big/ (q-1) \, \caluqMg \, $  is a
commutative  $ \C $--algebra  (cf.~[DP], \S 12).  Now from (3.1) we have
  $$  \runo = \prod_{\alpha \in R^+} {{\left( q_\alpha \cdot \ebar_\alpha
\otimes \fbar_\alpha;
\, q_\alpha^2 \right)}_\infty}^{\!\!\! -1} = \prod_{\alpha \in R^+}
\runoalpha  $$  letting  $ \, \runoalpha := {{\left( q_\alpha \cdot
\ebar_\alpha \otimes \fbar_\alpha; \, q_\alpha^2
\right)}_\infty}^{\!\!\! -1} \, $  for all  $ \, \alpha \in R^+ \, $,  and
Lemma 2.2 (for  $ \,
\varepsilon = 1 \, $) gives
  $$  \runoalpha = \exp\left({\,1\, \over \,q-1\,} \cdot {\,1\, \over \,2
d_\alpha\,} \cdot \varphi \left( q_\alpha \cdot \ebar_\alpha \otimes
\fbar_\alpha \right) \right) \cdot {\left( 1 - q_\alpha \cdot \ebar_\alpha
\otimes \fbar_\alpha \right)}^{-1/2} \cdot \big( 1 + {\Cal O}(q-1) \big)
$$  where we set  $ \, \varphi(z) := \sum_{n=1}^\infty {\,1\, \over
\,n^2\,} \, z^n \, $,  as usual.  Therefore we fall within the framework of
[Re], \S 3, hence we can apply Reshetikin's trick to conclude: namely, applying Lemma 3.2.2 of [Re] we get for all  $ \, \alpha \, \in R^+ \,$
  $$  \hbox{\rm Ad} \left( \runoalpha \right) (a) = \runoalpha \cdot a \cdot \runoalpha^{-1} \in
\caluqMg^{\widehat \otimes 2}  $$ for all  $ \, a \in \caluqMg^{\widehat \otimes 2} \, $,  i.~e.~$ \, \hbox{\rm Ad} \left( \runoalpha \right) $  restricts to an automorphism  $ \widetilde{\rcalunoalpha} $  of
$ \caluqMg^{\widehat \otimes 2} $;  thus also  $ \, \hbox{\rm Ad} \left( \runo \right) = \hbox{\rm Ad} \left( \prod_{\alpha \in R^+} \runoalpha \right) =
\prod_{\alpha \in R^+} \hbox{\rm Ad} \left( \runoalpha \right) = \prod_{\alpha \in R^+} \widetilde{\rcalunoalpha} \, $  does restrict to an automorphism  $ \widetilde{{\Cal  R}^{(1)}} $  of  $ \caluqMg^{\widehat \otimes 2} $  as claimed.  Then Theorem 3.2 yields the claim.   $ \square $
\enddemo

\vskip7pt

\proclaim{Corollary 3.8}  For any  $ c \in \C $,  let
$ {\Cal R}_c $  be the algebra automorphism  of  $ {{\Cal U}_c^M \! (\gerg)}^{\widehat \otimes 2} $  given by specialization of  $ \widetilde{\Cal R} $  at  $ q = c $.
Then  $ \left( {\Cal U}_c^M \! (\gerg), {\Cal R}_c \right) $  is a braided Hopf algebra.   $ \square $
\endproclaim
\eject

\vskip2,5truecm

\centerline {\bf \S \; 4 \  The geometrical meaning of the braiding structure
at roots of 1 }

\vskip10pt

   {\bf 4.1  Geometric framework.} \;  In this section we turn to geometry:
our aim is to show that the  series describing the adjoint action of the  $
R $-matrix  of a quantum group are more than formal objects, for they  {\sl
do converge},  in a proper sense, so that such action does yield
well-defined automorphisms of geometric objects.
                                                   \par
  Let  $ G $  be a connected simply connected semisimple Poisson algebraic
group over  $ \C $  with
$ \gerg $  as tangent Lie bialgebra; then there exists a uniquely defined
connected simply connected semisimple affine algebraic Poisson group  $ H $
over
$ \C $  with tangent Lie bialgebra  $ \gerg^* $  and algebra of polynomial
functions  $ F[H] $,  which is called  {\sl the Poisson group  {\sl dual}
of}  $ G $  (cf.~e.~g.~[DP], \S 11).
                                                  \par
   Let  $ \, \ell \in \N \, $  be  {\it odd},  $ \, \ell > d:= \max_i
\{d_i\} \, $,  or  $ \, \ell = 1
\, $;  then let  $ \, \varepsilon \in \C \, $  be a primitive  $ \ell $-th
root of 1.  As a matter of notation, let  $ \, U_\varepsilon :=
\caluepsilonPg \, $,  $ \, Z_\varepsilon := Z \left( U_\varepsilon
\right) \, $  (the centre of  $ U_\varepsilon $).  Everything in the sequel
can then be suitably extended to the case of quantum group  $ \uqMg $  with
general lattice  $ M $.  From the analysis in [DP] (cf.~also [DK], [DKP])
we recall the following results:
                                                    \par
   {\it (a)\/}  The subalgebra  $ \, Z_0 \, $  of
$ U_\varepsilon $  generated by  $ \, \ebar_\alpha^\ell $,  $
\fbar_\alpha^\ell $,  $ L_i^\ell \, $   ($ \alpha \in R^+ $,  $ i=\unon $)
is central, i.~e.~$ \, Z_0 \subseteq Z_\varepsilon \, $.
                                                       \par
   {\it (b)\/}  $ Z_\varepsilon $  and  $ Z_0 $  inherit (from  $ \caluqPg
$)  canonical structures of  Poisson algebras; in particular,  $ Z_0 $  is
a Poisson Hopf algebra.
                                                       \par
   {\it (c)\/}  There exists an isomorphism  $ \, Z_0 \cong F[H] \, $  as
Poisson Hopf algebras (with  respect to a suitable normalization of the
Poisson bracket on  $ Z_ 0 $), hence  $ \, Spec \big( Z_0
\big) \cong H \, $  as Poisson (complex affine algebraic) groups.  In
particular (for  $ \, \ell = 1
\, $)  $ \, Spec \left( \caluunoPg \right) \cong H \, $.
                                                        \par
   Recall that in [DK], [DKP], [DP] the spectra  $ \, Spec \big(
U_\varepsilon \big) $,  $ Spec \big(  Z_\varepsilon \big) $,  and  $ Spec
\big( Z_0 \big) $  are introduced as the set of isomorphism classes of
finite dimensional representations of the corresponding algebras  $
U_\varepsilon $,  $ Z_\varepsilon
$,  and  $ Z_0 $;  in particular  $ Spec \big( Z_\varepsilon \big) $  and
$ Spec \big( Z_0 \big) $  can be identified with usual geometric objects,
namely complex affine algebraic varieties describing the maximal spectrum
of  $ Z_\varepsilon $  and  $ Z_0 $;  since  $ \, Z_0 \cong F[H] \, $  as
Poisson Hopf algebras, we also have  $ \, Spec \big( Z_0 \big) \cong H \,
$  as Poisson affine algebraic groups (over  $ \C $);  thus in the sequel
we will also set  $ \, H_\varepsilon := Spec \big( Z_\varepsilon
\big) \, $  and  $ \, S_\varepsilon := Spec \big( U_\varepsilon \big) \,
$.  The analysis in [DP] describes  $ Spec \big( U_\varepsilon \big) $  as
({\it espace \'etal\'e\/}  of) a sheaf   --- or a fibre bundle ---   of
algebras over  $ Spec \big( Z_0 \big) $  or  $ Spec \big( Z_\varepsilon
\big) $;  in particular we can think at  $ U_\varepsilon $  as the algebra
of global sections of this sheaf.
                                                       \par
   Now set
  $$  y_\alpha := \fbar_\alpha^\ell {\Big\vert}_{q=\varepsilon} \, ,
\qquad  z_\lambda :=  L_\lambda^\ell {\Big\vert}_{q=\varepsilon} \, ,
\qquad x_\alpha := \ebar_\alpha^\ell {\Big\vert}_{q=\varepsilon}   \;
\qquad  \forall\, \alpha \in R^+, \, \lambda \in P  $$  and in particular
$ \, y_i := y_{\alpha_i} \, $,  $ \, z_i := z_{\omega_i} \, $,  $ \, x_i :=
x_{\alpha_i} \, $  ($ i=\unon $).  Following [DK], \S 3.5, we denote by  $
\widehat{Z}_0 $  the algebra of all formal power series in the  $ y_\alpha
$'s,  $ z_i^{\pm 1} $'s,  $ x_\alpha $'s  which converge to meromorphic
functions for all complex values of the  $ y_\alpha $'s,  $ x_\alpha $'s,
and all non-zero complex values of the  $ z_i $'s;  then let  $ \,
\widehat{U}_\varepsilon := \widehat{Z}_0
\otimes_{Z_0} U_\varepsilon $,  $ \widehat{Z}_\varepsilon := \widehat{Z}_0
\otimes_{Z_0} Z_\varepsilon
\, $.  In other words we can think at  $ \widehat{U}_\varepsilon $  as the
algebra of global meromorphic sections of the corresponding bundle of
algebras over  $ \, Spec \big( Z_0 \big) \cong H \, $.  Similar notations
and definitions will be used when dealing with square tensor powers, like
$ {Z_0}^{\otimes 2}
$,  $ {Z_\varepsilon}^{\otimes 2} $,  and so on.  Notice also that  $ \,
Spec \big( Z_0 \otimes Z_0
\big) = Spec \big( Z_0 \big) \times Spec \big( Z_0 \big) = H \times H \,
$,  $ \, Spec \big( Z_\varepsilon \otimes Z_\varepsilon \big) = Spec \big(
Z_\varepsilon \big) \times Spec \big( Z_\varepsilon \big) = H_\varepsilon
\times H_\varepsilon \, $,  $ \, Spec(U_\varepsilon \otimes U_\varepsilon)
= Spec(U_\varepsilon) \times Spec(U_\varepsilon) = S_\varepsilon \times
S_\varepsilon
\, $,
                                       \par
  {\sl Warning:\/}  when dealing with cross-product spaces like  $ X \times
Y $,  we shall use left  subscripts to denote functions of either of the
two spaces, viz.~$ \, {_2}x := 1 \otimes x \, $,  $ \, {_1\ebar_\alpha} :=
\ebar_\alpha \otimes 1 \, $,  etc.
                                              \par
   Let  $ {\Cal H}^{(N)} $  be any fixed ramified  $ N $-fold  covering
(for  $ \, N \in \N \cup
\{ \infty \} \, $)  of  $ H $  (so that  $ {\Cal H}^{(N)} \times {\Cal
H}^{(N)} $  is an  $ N $-fold  covering of  $ H \times H \, $);  then we
denote by  $ {\Cal H}_\varepsilon^{(N)} $  and  $ {\Cal
S}_\varepsilon^{(N)} $  the fiber products  $ \, {\Cal H}_\varepsilon^{(N)}
:= {\Cal H}^{(N)} \times_H H_\varepsilon \, $  and  $ \, {\Cal
S}_\varepsilon^{(N)} := {\Cal H}^{(N)} \times_H S_\varepsilon \, $.  Notice
that  $ \, {\Cal H}_\varepsilon^{(N)} = Spec \left( \widehat{Z}_\varepsilon
\right) \, $  and
$ \, {\Cal S}_\varepsilon^{(N)} = Spec \left( \widehat{U}_\varepsilon
\right) \, $;  furthermore,
$ {\Cal H}^{(N)} $  and  $ {\Cal H}_\varepsilon^{(N)} $  clearly have a
unique Poisson structure compatible with the covering map, so that  $ {\Cal
H}_\varepsilon^{(N)} $  is a (complex analytic) Poisson variety and  $
{\Cal H}^{(N)} $  is a (complex analytic) Poisson group.  Finally,  $ \,
\tau :=
\sigma^* \colon H \times H \longrightarrow H \times H \, $  ($ \sigma $
being defined in \S 1.3) is given by  $ \, (x,y) \mapsto (y,x) \, $;  then
$ \, \tau^{(N)} \colon {\Cal H}^{(N)} \times {\Cal H}^{(N)} \longrightarrow
{\Cal H}^{(N)} \times {\Cal H}^{(N)} \, $,  also given by  $ \, (x,y)
\mapsto (y,x) \, $,  is a lifting of  $ \tau $  to  $ {\Cal H}^{(N)} \times
{\Cal H}^{(N)} $.
                                              \par
   Fix now  $ \, \ell > 1 \, $:  we are ready for the next result, which
claims that the "formal  automorphism"  $ {\Cal R}_\varepsilon $  giving
the braiding structure of  $ U_\varepsilon $  actually does
converge in a proper sense.

\vskip7pt

\proclaim {Proposition 4.2}  The algebra automorphism  $ \, {\Cal R}_\varepsilon \colon  U_\varepsilon^{\widehat \otimes 2} @>>> U_\varepsilon^{\widehat \otimes 2} \, $  defines a meromorphic automorphism  $ {\Cal R}_{\varepsilon,\infty}^\ast $  of  $ \,  {\Cal
S}_\varepsilon^{(\infty)} \times {\Cal S}_\varepsilon^{(\infty)} \, $,  which restricts to meromorphic Poisson automorphisms  $ {\Cal
R}_{\varepsilon,\infty}^\ast \colon \, {\Cal H}_\varepsilon^{(\infty)} \times {\Cal H}_\varepsilon^{(\infty)} \longrightarrow {\Cal
H}_\varepsilon^{(\infty)} \times {\Cal H}_\varepsilon^{(\infty)} \, $  and  $ {\Cal R}_{\varepsilon,\infty}^\ast \colon \, {\Cal H}^{(\infty)} \times {\Cal H}^{(\infty)} \longrightarrow {\Cal H}^{(\infty)} \times {\Cal H}^{(\infty)} \, $.
                                           \hfill\break
   \indent   Moreover,  $ {\Cal R}_{\varepsilon,
\infty}^\ast $  and its restrictions enjoy the dual properties of (1.4--6);
in particular,
$ \, {\Cal R}_{\varepsilon,\infty}^\ast \neq \tau \, $,  and  $ \, m \big(
{\Cal R}_{\Cal H}^\ast (x,y)
\big) = m(y,x) = y \cdot x \, $  for all  $ \, x, y \in {\Cal
H}^{(\infty)}
\, $  ($ m $  and  "$ \, \cdot \, $"  denoting the product of
$ {\Cal H}^{(\infty)} \, $),  and a braid group action exists on
$ \times $--powers  of  $ {\Cal H}^{(\infty)} $.
\endproclaim

\demo{Proof}  The first step in the proof amounts to show that series  $ \,
{\Cal R}_\varepsilon (x \otimes y) \, $  do converge almost everywhere on a
suitable covering  $ \, {\Cal S}_\varepsilon^{(\infty)} \times {\Cal
S}_\varepsilon^{(\infty)} \, $.  Recall (cf.~Proposition~3.2 and its proof)
that  $ \, \widetilde{\Cal R} :=
\widetilde\rcaluno \smallcirc \widetilde\rcalzero \, $,  thus  $ \, {\Cal
R}_\varepsilon := \rcaluno_{\!\! \varepsilon} \smallcirc \rcalzero_{\!\!
\varepsilon} \, $  with   $ \,
\rcalzero_{\!\! \varepsilon} := \widetilde\rcalzero \mod (q-\varepsilon) \,
$  and  $ \,
\rcaluno_{\!\! \varepsilon} := \widetilde\rcaluno \mod (q-\varepsilon) \,
$.  For  $ \rcalzero_{\!\! \varepsilon} $  the very definition implies that
no problem of convergence (nor of domain of definition) occurs.  For  $
\rcaluno_{\!\! \varepsilon} $,  recall that  $ \, \widetilde\rcaluno :=
\hbox{\rm Ad} \left( \runo \right) {\Big\vert}_{\caluqPg^{\widehat \otimes
2}} \, $  and
  $$  \runo := \prod_{\alpha \in R^+} {{\left( q_\alpha \cdot \ebar_\alpha
\otimes \fbar_\alpha; \, q_\alpha^2 \right)}_\infty}^{\!\!\! -1} =
\prod_{\alpha \in  R^+} \runoalpha   \eqno (4.1)  $$ where  $ \, \runoalpha
:= {{\left (q_\alpha \cdot \ebar_\alpha \otimes \fbar_\alpha; \, q_\alpha^2
\right)}_\infty}^{\!\!\! -1} \, $,  like in the proof of Proposition~3.2;
therefore
  $$  \hbox{\rm Ad} \left( \runo \right) = \prod_{\alpha \in R^+} \hbox{\rm
Ad} \left( \runoalpha \right) = \prod_{\alpha \in R^+} \hbox{\rm Ad} \left(
{{\left( q_\alpha \cdot \ebar_\alpha \otimes \fbar_\alpha; \, q_\alpha^2
\right)}_\infty}^{\!\!\! -1} \right) \, .   \eqno (4.2)  $$
   \indent   Now again we apply Reshetikin's trick: from [Re], Lemma 3.2.2
and formulas (3.2.10--11), and from our Lemma~2.2 we get
  $$  \eqalign{
   \hbox{\rm Ad} \left(\runoalpha \right) \mod (q-\varepsilon) =  &
\hbox{\rm Ad} \left( {{\left( q_\alpha \cdot \ebar_\alpha \otimes
\fbar_\alpha; \, q_\alpha^2 \right)}_\infty}^{\!\!\! -1} \right) \mod
(q-\varepsilon) =  \cr
   = \hbox{\rm Ad}\left( \exp\left({\Phi_\alpha \over q-\varepsilon}\right)
\right)  &  \mod (q-\varepsilon) =
\exp \Big( \hbox{\rm ad}_{\{ \ , \ \}} \big( \Phi_\alpha \big) \Big) \mod
(q-\varepsilon) =  \cr
   = \exp \bigg( \hbox{\rm ad}_{\{ \ , \ \}} \bigg( {\,\varepsilon\, \over
\,2 d_\alpha \ell^{\scriptstyle 2}\,} \, \cdot  &  \, \varphi \big(
q_\alpha \ebar_\alpha \otimes \fbar_\alpha \big) \bigg) \bigg) \mod
(q-\varepsilon)  \cr }
\eqno (4.3)  $$
whith
  $$  \eqalign{
   \Phi_\alpha  &  := \left( {q-\varepsilon \over
q_\alpha^{2 \ell^{\scriptscriptstyle 2}} - 1} \cdot \varphi \left(
q_\alpha \ebar_\alpha \otimes \fbar_\alpha \right) \, - (q-\varepsilon)
\cdot \log \left( \prod_{k=0}^{\ell - 1} {\,k\, \over \,\ell\,} {\left( 1 -
\varepsilon^k q_\alpha \ebar_\alpha \otimes \fbar_\alpha \right)} \right)
\right) \cdot  \cr
                &  \cdot \big( 1 + {\Cal O}(q-\varepsilon) \big) = \exp
\left( {q-\varepsilon \over q_\alpha^{2 \ell^{\scriptscriptstyle 2}} - 1} \cdot \varphi \left( q_\alpha
\ebar_\alpha \otimes \fbar_\alpha \right) \right)  \qquad  \mod (q -
\varepsilon)  \cr }  $$
and
  $$  \varphi(t) := \int_0^{t^\ell} {\log(1-\tau) \over \tau} \, d\tau =
\sum_{n=1}^\infty {t^{\ell n} \over n^2} \qquad \hbox{ (by Taylor
expansion) } \, .  $$
   \indent   Notice that
  $$  \displaylines{
   \hbox{ad}_{\{ \ , \ \}} \big( t \cdot \varphi(x) \big) (y) = \big\{ t
\cdot \varphi(x), y \big\} = t \cdot \left\{ \sum_{n=1}^\infty {x^{\ell n}
\over n^2}, y \right\} = t \cdot \sum_{n=1}^\infty {1 \over n^2} \cdot
\left\{ {\left( x^\ell \right)}^n, y \right\} = \hfill \quad  \cr
   {} \quad \hfill  = t \cdot \sum_{n=1}^\infty {\,1\, \over \,n^2\,} \, n
{\left( x^\ell \right)}^{n-1} \cdot \left\{ x^\ell, y \right\} =
\sum_{n=1}^\infty {{\left( x^\ell \right)}^{n-1} \over n} \cdot \left\{ t
\cdot x^\ell, y \right\}  \cr }  $$
(because of Leibnitz' rule:  $ \, \{ \,
\cdot \, , y\} = -
\hbox{ad}_{\{ \ , \ \}}(y)  \, $  is a derivation!), hence
  $$  \hbox{ad}_{\{ \ , \ \}} \big( t \cdot \varphi(x) \big) = \psi \left(
x^\ell \right) \cdot \hbox{ad}_{\{ \ , \ \}} \left( t \cdot x^\ell \right)
$$ with  $ \, \psi(t):= {\,\log(1-y)\, \over \,y\,} = \sum_{n=0}^\infty
{y^n \over n+1} \, $  (by Taylor expansion again), and then
  $$  \exp \! \Big( \hbox{ad}_{\{ \ , \ \}} \big( t \cdot \varphi(x) \big)
\Big) = \exp \! \Big( \psi \big( x^\ell \big) \cdot \hbox{ad}_{\{ \ , \ \}}
\big( t \cdot x^\ell \big) \Big) \; ;  $$ together with (4.3) this gives
  $$  \displaylines{
   \hbox{Ad} \left( \runoalpha \right) \mod (q-\varepsilon) = \exp\left(
\hbox{ad}_{\{ \ , \ \}} \! \left( {\,\varepsilon\, \over \,2 d_\alpha
\ell^{\scriptstyle 2}\,} \cdot \varphi \left( q_\alpha \ebar_\alpha \otimes
\fbar_\alpha \right) \right) \right) =   \hfill \quad {}  \cr
  = \exp \left( \psi \left( q_\alpha^\ell \ebar_\alpha^\ell \otimes
\fbar_\alpha^\ell \right) \cdot \hbox{ad}_{\{ \ , \ \}} \! \left(
\,{\varepsilon\, \over \,2 d_\alpha \ell^{\scriptstyle 2}\,} \,
q_\alpha^\ell \cdot
\ebar_\alpha^\ell \otimes  \fbar_\alpha^\ell \right) \right) \mod
(q-\varepsilon) =  \cr
   {} \quad \hfill   = \exp \left( {\,\varepsilon\, \over \,2 d_\alpha
\ell^{\scriptstyle 2}\,} \, \psi \!\left( x_\alpha \otimes y_\alpha \right)
\cdot \hbox{ad}_{\{ \ , \ \}} \!\left( x_\alpha \otimes y_\alpha \right)
\right) \mod (q- \varepsilon) \, .  \cr }  $$
   \indent   Therefore we have to show that the formal series
  $$  \exp \left( {\,\varepsilon\, \over \,2 d_\alpha
\ell^{\scriptscriptstyle 2}\,} \cdot \psi \!\left( x_\alpha \otimes
y_\alpha \right) \cdot \hbox{ad}_{\{ \ , \ \}} \!\left( x_\alpha \otimes
y_\alpha \right) \right) (x \otimes y)  $$  for  $ x $,  $ y $  generators
of  $ U_\varepsilon $  (that is  $ \, x, y \in \big\{\, 1, \, \fbar_\alpha,
\, L_i, \, \ebar_\alpha  \mod (q-\varepsilon) \,\big\vert\, i=\unon, \alpha
\in R^+ \,\big\} \, $)  does converge (to a meromorphic function on  $
{\Cal S}_\varepsilon^{(\infty)} \times {\Cal S}_\varepsilon^{(\infty)} \,
$).  But notice that the following obvious identity holds (for all  $ \, n
\in \N \, $)
  $$  {\left( {\,\varepsilon\, \over \,2 d_\alpha \ell^{\scriptscriptstyle
2}\,} \cdot
\psi \!\left( x_\alpha \otimes y_\alpha \right) \cdot \hbox{ad}_{\{ \ , \
\}} \!\left( x_\alpha \otimes y_\alpha \right) \right)}^n = \, {\psi
\!\left( x_\alpha \otimes y_\alpha \right)}^n \cdot {\left(
{\,\varepsilon\, \over \,2 d_\alpha \ell^{\scriptscriptstyle 2}\,} \cdot
\hbox{ad}_{\{ \ , \ \}}
\!\left( x_\alpha \otimes y_\alpha \right) \right)}^n  $$  because of
Leibnitz' rule and  $ \, \hbox{\rm ad}_{\{ \ , \ \}} \!\left( x_\alpha
\otimes y_\alpha
\right) \left( {\,\varepsilon\, \over \,2 d_\alpha \ell^{\scriptscriptstyle
2}\,} \, \psi \!\left( x_\alpha \otimes y_\alpha \right)
                         \right) = 0 \, $;  moreover,\break
$ \, \psi \!\left( x_\alpha \otimes y_\alpha \right) \, $  is a meromorphic
function on the
$ \infty $-fold  ramified covering  $ {\Cal H}^{(\infty)} \times {\Cal
H}^{(\infty)} $  of  $ H \times  H $.  Now recall that
  $$  [x \otimes y, z \otimes w] = [x, z]  \otimes y w + x z \otimes [y, w]
\; ;   \eqno (4.4)  $$
then set  $ \, {_1\hbox{e}_\alpha} := \hbox{ad}_{[ \
, \  ]} \left( {_1E_\alpha^{(\ell)}}
\right){\Big\vert}_{q=\varepsilon} \, $,  $ \, {_2\hbox{f}_\alpha} := \hbox{ad}_{[ \ , \ ]} \left( {_2F_\alpha^{(\ell)}} \right){\Big\vert}_{q=\varepsilon}
\, $  ($ \, X_\alpha^{(n)} := { X_\alpha^n \over {[n]}_{q_\alpha}! } \, $),  observe that
  $$  \displaylines{
   {_1\hbox{e}_\alpha} := \hbox{ad}_{[ \ , \  ]} \left(
{_1E_\alpha^{(\ell)}} \right)  {\Big\vert}_{q=\varepsilon} = {\left( {\,1\,
\over \, q_\alpha^{2 \ell^{\scriptscriptstyle 2}} - 1 \,}
\cdot \hbox{ad}_{[ \ , \  ]} \left( {_1\ebar_\alpha^\ell} \right)
\right)}{\bigg\vert}_{q=\varepsilon} = {\,\varepsilon\, \over \,2 d_\alpha
{\ell}^{\scriptstyle 2}\,}
\cdot \hbox{ad}_{\{ \ , \ \}} \left( {_1x_\alpha} \right)  \cr
   {_2\hbox{f}_\alpha} := \hbox{ad}_{[ \ , \  ]} \left(
{_2F_\alpha^{(\ell)}} \right) {\Big\vert}_{q=\varepsilon} = {\left( {\,1\,
\over \, q_\alpha^{2 {\ell}^{\scriptscriptstyle 2}} - 1 \,} \cdot
\hbox{ad}_{[ \ , \  ]} \left( {_2\fbar_\alpha^\ell} \right)
\right)}{\bigg\vert}_{q=\varepsilon} = {\,\varepsilon\, \over \,2 d_\alpha
\ell^{\scriptstyle 2}\,} \cdot \hbox{ad}_{\{ \ , \ \}} \left( {_2y_\alpha}
\right)  \cr }  $$  and let  $ \, m(x) \colon y \mapsto x y \, $  (left
multiplication by  $ x $);  then formula (4.4) gives
  $$  {\,\varepsilon\, \over \,2 d_\alpha \ell^{\scriptscriptstyle 2}\,}
\cdot \hbox{ad}_{\{ \ , \ \}}
\left( x_\alpha \otimes  y_\alpha \right) = {_1\hbox{e}_\alpha} \otimes m
\left( {_2y_\alpha} \right) + m \left( {_1x_\alpha} \right) \otimes
{_2\hbox{f}_\alpha} \, ;   \eqno (4.5)  $$   one trivially checks that  $
\, {_1\hbox{e}_\alpha} \otimes m \left( {_2y_\alpha} \right) \, $  and
$ \, m \left( {_1x_\alpha} \right) \otimes {_2\hbox{f}_\alpha} \, $  are
operators which commute with each other, thus (4.5) gives
  $$  \exp \! \Big( \hbox{ad}_{\{ \ , \ \}} \! \left( x_\alpha \otimes
y_\alpha \right) \Big) =
\exp \! \big( {_1\hbox{e}_\alpha} \otimes m \left( {_2y_\alpha} \right)
\big) \smallcirc \exp \! \big( m
\left( {_1x_\alpha} \right) \otimes {_2\hbox{f}_\alpha} \big) \, ;   \eqno
(4.6)  $$  for  $ x $,  $ y $  generators of  $ U_\varepsilon $  we have
  $$  \eqalign{
   &  \exp \! \big( \psi \! \left(\! x_\alpha \!\otimes\! y_\alpha
\!\right) \cdot  m \left(\! x_\alpha \!\right) \!\otimes\! \hbox{f}_\alpha
\big) (x \!\otimes\! y) = \exp \! \left( { \,\log
\left(\! 1 - {_1x_\alpha} \!\cdot\! {_2y_\alpha} \!\right)\, \over
\,{_1x_\alpha} \cdot {_2y_\alpha}\, } \cdot {_1x_\alpha} \cdot
\hbox{f}_\alpha \right) \! (x \!\otimes\! y) =  \cr
   &  = \exp \! \left( { \,\log \left( 1 - {_1x_\alpha} \cdot {_2y_\alpha}
\right)\, \over
\,{_2y_\alpha}\, } \cdot \hbox{f}_\alpha \right) \! (x \!\otimes\! y) =
\exp \! \left( { \,\log \left( 1 -  {_1x_\alpha} \cdot {_2y_\alpha}
\right)\, \over \,{_2y_\alpha}\, } \cdot \hbox{f}_\alpha \right) \! (y)
\cdot x \! \otimes \! 1 .  \cr }   \eqno (4.7)  $$
   \indent   It is proved in [DK], \S 3, that  $ \exp \left( t \cdot
\hbox{f}_\alpha \right) $   converges to a holomorphic automorphism of the
algebra of global holomorphic sections of
$ S_\varepsilon $  (as a bundle over  $ H $),  for all  $ \, t \in \C \,
$;  when  $ t $  is replaced with any meromorphic function on  $ H $,  the
series we get does converge to an automorphism of the algebra of  {\sl
meromorphic}  sections (cf.~formulas in the proof of Proposition 3.5 of
[DK]); since
$ \, { \,\log \left( 1 - {_1x_\alpha} \cdot {_2y_\alpha} \right)\, \over
\,{_2y_\alpha}\, } \, $  is meromorphic on the  $ \infty $-fold  covering
$ {\Cal H}^{(\infty)} \times {\Cal H}^{(\infty)} $,  we conclude that  $ \,
\exp \left( \psi \!\left( x_\alpha \otimes y_\alpha \right) \cdot m\!\left(
{_1x_\alpha} \right) \otimes {_2\hbox{f}_\alpha} \right) \! (x \!\otimes\! y)
\, $  is a meromorphic section of  $ \, {\Cal S}_\varepsilon^{(\infty)}
\times {\Cal S}_\varepsilon^{(\infty)} \, $;  the same holds for  $ \, \exp
\left( \psi \!\left( x_\alpha \otimes y_\alpha \right) \cdot
{_1\hbox{e}_\alpha} \otimes m\!\left( {_2y_\alpha} \right) \right) \, $,
and finally for  $ \, \hbox{Ad} \left( \runoalpha \right)
{\Big\vert}_{q=\varepsilon} = \exp \left( {\,\varepsilon\, \over \,2
d_\alpha \ell^{\scriptscriptstyle 2}\,} \cdot \psi \!\left( x_\alpha
\otimes y_\alpha \right) \cdot \hbox{ad}_{\{\ ,\ \}} \!\left( x_\alpha
\otimes y_\alpha \right) \right) {\Big\vert}_{q=\varepsilon} \, $,  q.~e.~d.
                                          \par
   For the second part, notice that  $ \rcalzero_{\!\!\! \varepsilon} $
clearly leaves invariant both
$ \widehat{ {Z_\varepsilon}^{\otimes 2}} $  and  $ \widehat{{Z_0}^{\otimes
2}} $,  hence its dual leaves invariant  $ \, {\Cal
H}_\varepsilon^{(\infty)} \times {\Cal H}_\varepsilon^{(\infty)} \, $  and
$ \, {\Cal H}^{(\infty)} \times {\Cal H}^{(\infty)} \, $;  moreover, since
$ \rcaluno_{\!\! \varepsilon} $  is a product of terms
  $$  \hbox{Ad} \left( \runoalpha \right) {\Big\vert}_{q=\varepsilon} =
\exp \left( {\,-\varepsilon\,
\over \,2 d_\alpha \ell^{\scriptstyle 2}\,} \cdot \psi \left( x_\alpha
\otimes y_\alpha \right) \cdot
\hbox{ad}_{\{ \ , \ \}} \left( x_\alpha \otimes y_\alpha \right) \right)
{\bigg\vert}_{q=\varepsilon}  $$ since  $ \widehat{{Z_\varepsilon}^{\otimes
2}} $  and  $ \widehat{{Z_0}^{\otimes 2}} $  are closed for the Poisson
bracket, and since  $ \, x_\alpha \otimes y_\alpha \in
\widehat{{Z_0}^{\otimes 2}} \subseteq
\widehat{{Z_\varepsilon}^{\otimes 2}} \, $,  we have that the dual of
$ \rcaluno_{\!\!\! \varepsilon} $  leaves  $ \, {\Cal
H}_\varepsilon^{(\infty)} \times  {\Cal H}_\varepsilon^{(\infty)} \, $
and  $ \, {\Cal H}^{(\infty)} \times {\Cal H}^{(\infty)} \, $  invariant;
thus we conclude that  $ {\Cal R}_\varepsilon^\ast $  leaves  $ \, {\Cal
H}_\varepsilon^{(\infty)} \times {\Cal H}_\varepsilon^{(\infty)} \, $  and
$ \, {\Cal H}^{(\infty)} \times {\Cal H}^{(\infty)} \, $  invariant.
Finally, it clearly preserves the Poisson structure because  $ {\Cal
R}_\varepsilon $  is defined by specializing an algebra automorphism of  $
\widehat{\caluqPg^{\otimes 2}} $,  whence
  $$  {\Cal R}_\varepsilon \big( \{x_0, y_0\} \big) = \widetilde{\Cal
R}\left( [x,y] \over  q-\varepsilon \right) \! {\Bigg\vert}_{q=\varepsilon}
= {\,\left[ \widetilde{\Cal R}(x), \widetilde{\Cal R}(y) \right]\, \over
\,q-\varepsilon\,} {\Bigg\vert}_{q=\varepsilon} = \big\{ {\Cal
R}_\varepsilon(x_0), {\Cal R}_\varepsilon(y_0) \big\} \, . $$
   \indent   The proof of the last part of the statement is completely
trivial, by functoriality.   $ \square $
\enddemo

\vskip7pt

   A deeper analysis yelds to improve the previous result, proving that the
convergence already holds  on  {\sl finite}  ramified coverings, as the
following shows.

\vskip7pt

\proclaim{Theorem 4.3}  The meromorphic automorphism  $ \, {\Cal
R}_\varepsilon^\ast \colon  {\Cal S}_\varepsilon^{(\infty)} \times {\Cal
S}_\varepsilon^{(\infty)} \longrightarrow {\Cal S}_\varepsilon^{(\infty)}
\times {\Cal S}_\varepsilon^{(\infty)} \, $  pushes down to a birational
automorphism  $ \, {\Cal R}_{\varepsilon, \ell}^{\,\ast} \colon {\Cal
S}_\varepsilon^{(2 \ell)} \times {\Cal S}_\varepsilon^{(2 \ell)}
\longrightarrow {\Cal S}_\varepsilon^{(2 \ell)} \times {\Cal
S}_\varepsilon^{(2 \ell)} \, $;  moreover,  $ \, {\Cal R}_{\varepsilon,
\ell}^{\,\ast} \neq \tau^{(2\ell)} $,  and  $ {\Cal R}_{\varepsilon, \ell}^{\,\ast} $  enjoys the dual
properties of (1.4--6).
                                     \hfill\break
   \indent   The same holds with  $ H_\varepsilon $,  resp.~$ H $  instead
of  $ S_\varepsilon $,  with  a birational Poisson automorphism  $ \, {\Cal
R}_{\varepsilon, \ell}^{\,\ast} \colon {\Cal H}_\varepsilon^{(2 \ell)}
\times {\Cal H}_\varepsilon^{(2 \ell)} \rightarrow {\Cal H}_\varepsilon^{(2
\ell)} \times {\Cal H}_\varepsilon^{(2 \ell)} \, $,  resp.~$ \, {\Cal
R}_{\varepsilon, \ell}^{\,\ast}
\colon {\Cal H}^{(2 \ell)} \times {\Cal H}^{(2 \ell)} \rightarrow {\Cal
H}^{(2 \ell)} \times {\Cal H}^{(2 \ell)} \, $:  in particular,  $ \, m
\left( {\Cal R}_{\varepsilon, \ell}^{\,\ast} (x,y) \right) = m(y,x) = y
\cdot x \, $  for all  $ \, x, y \in {\Cal H}^{(2 \ell)} \, $  (where  $ m $
and  "$ \, \cdot \, $"  denote the product of  $ \, {\Cal H}^{(2 \ell)} \, $),
and a braid group action exists on  $ \times $--powers  of
$ {\Cal H}^{(2 \ell)} $.
\endproclaim

\demo{Proof}  It is clear that for  $ \rcalzero_{\!\!\! \varepsilon} $
everything is o.k.  As for  $ \rcaluno_{\!\!\! \varepsilon} $,  from the
proof of Proposition 4.2 we see that it is enough to show that
  $$  \exp \left( {\,- \varepsilon\, \over \,2 d_\alpha
\ell^{\scriptscriptstyle 2}\,}
\psi \!\left( x_\alpha \otimes y_\alpha \right) \cdot \hbox{ad}_{\{\ ,\ \}}
\left( x_\alpha \otimes y_\alpha \right) \right) (x \otimes y)   \eqno
(4.8)  $$  (for any  $ x $,  $ y $  in  $ U_\varepsilon $)  is a rational
section of the bundle  $ {\Cal S}_\varepsilon^{(2 \ell)} $  on a  $ 2 \ell
$-fold  ramified covering  $ {\Cal H}^{(2 \ell)} $  of  $ H \, $:  this
again amounts to perform some computations.  In particular (cf.~(4.4--7))
we are reduced to check the same for functions
  $$  \eqalign{
   \exp  &  \left( {\,\log \left( 1 - {_1x_\alpha} \cdot {_2y_\alpha}
\right)\, \over \,{_1x_\alpha}\,} \cdot {_1\hbox{e}_\alpha}\, \right) \big(
{_1x} \cdot {_2y} \big)  \cr
   \exp  &  \left( {\,\log \left( 1 - {_1x_\alpha} \cdot {_2y_\alpha}
\right)\, \over \,{_2y_\alpha}\,} \cdot {_2\hbox{f}_\alpha}\, \right) \big(
{_1x} \cdot {_2y} \big)  \cr }   \eqno (4.9)  $$
   \indent   We deal with the first function above, the proof for the
second following by symmetry.
                                                 \par
   Since  $ \, {\,\log \left( 1 - {_1x_\alpha} \cdot {_2y_\alpha} \right)\,
\over \,{_1x_\alpha}\,} \cdot {_1\hbox{e}_\alpha} \, $  is a derivation of
$ \widehat{{U_\varepsilon}^{\otimes 2}} $,  its exponential is an
automorphism of  $ \widehat{{U_\varepsilon}^{\otimes 2}} \, $;  now  $ \,
{\,\log \left( 1 - {_1x_\alpha} \cdot {_2y_\alpha} \right)\, \over
\,{_1x_\alpha}\,} \cdot {_1\hbox{e}_\alpha} (1 \otimes y) = 0 \, $,
whence  $ \, \exp \left( {\,\log \left( 1 - {_1x_\alpha} \cdot {_2y_\alpha}
\right)\, \over \,{_1x_\alpha}\,} \cdot {_1\hbox{e}_\alpha} \right) (1
\otimes y) =  1 \otimes y
\, $,  for all  $ \, y \in U_\varepsilon $;  therefore we have only to
compute  $ \, \exp \left( {\,\log \left( 1 - {_1x_\alpha} \cdot
{_2y_\alpha} \right)\, \over \,{_1x_\alpha}\,} \cdot {_1\hbox{e}_\alpha}
\right) (x \otimes 1) \, $  for  $ \, x = {_1x} \in U_\varepsilon \, $:  in
particular, it is enough to take  $ x $  to be a generator of  $
U_\varepsilon $,  namely  $ \, x \in \big\{\, F_i, L_\lambda, E_j
\,\big\vert\, i, j= \unon; \lambda \in P \,\big\} \, $.
                                                       \par
   Like in the proof of [DK], Proposition 3.5, exploiting the braid group
action we can restrict to the case of simple roots  $ \, \alpha = \alpha_i
\, $,  $ \, i = \unon \, $  (thus we set  $ \, {_1\ebar_i} :=
{_1\ebar_{\alpha_i}} $,
$ {_1\hbox{e}_i} := {_1\hbox{e}_{\alpha_i}} $,  and so on),  using formulas
   $$
    \matrix
   {_1\ebar_\alpha} = T_w \left( {_1\ebar_i} \right) \, ,  &  \qquad
{_1\fbar_\alpha} =  T_w \left( {_1\fbar_i} \right)  &  \qquad \qquad
\hbox{ for \ } \, \alpha = w(\alpha_i)  \\
   {_1\hbox{e}_\alpha} = T_w \smallcirc {_1\hbox{e}_i} \smallcirc T_w^{-1}
\, ,  &  \qquad  {_1\hbox{f}_\alpha} =  T_w \smallcirc {_1\hbox{f}_i}
\smallcirc T_w^{-1}  &  \qquad \qquad  \hbox{ for
\ } \, \alpha = w(\alpha_i)  \\
    \endmatrix   \eqno (4.10)  $$ (cf.~[DK], \S 3.4), where  $ T_w $
denotes the unique element of the braid group associated to
$ w \in W $.  Moreover, from direct computation or resuming formulas in the
proof of [DK], Proposition 3.5, we get,  {\it mutatis mutandis},
  $$  \eqalign{
  \exp  &  \left( t \cdot {_1\hbox{e}_i} \right) \left( {_1L_\lambda}
\right) = e^{ (\langle \alpha_i \vert \lambda
\rangle / 2 \ell) \cdot t \cdot {_1\!x_i}} \cdot {_1L_\lambda}  \cr
  \exp  & \left( t \cdot {_1\hbox{e}_i} \right) \left( {_1\fbar_j} \right)
=  \cr
        &  = {_1\fbar_j} - \delta_{i j} \left( {\,{e^{-t \cdot {_1\!x_i} /
\ell} - 1}\, \over \,{_1x_i}\,} \cdot \varepsilon^{d_i} L_{\alpha_i} +
{\,{e^{t \cdot {_1\!x_i} / \ell} - 1}\, \over \,{_1x_i}\,} \cdot
\varepsilon^{-d_i} L_{-\alpha_i} \right) \cdot {_1\ebar_i^{\ell - 1}}
\cr }  $$  for any indeterminate  $ t $  which commutes with  $ {_1\ebar_i}
$  (where  $ \, \langle \alpha_i \vert
\lambda \rangle := 2 ( \alpha_i \vert \lambda) /  (\alpha_i \vert \alpha_i)
\in \Z \, $);  when instead  of  $ t $  we have the meromomorphic function
$ \, {\,\log \left( 1 - {_1x_i} \cdot {_2y_i} \right)\,
\over \,{_1x_i}\,} \, $  (which does commute with  $ {_1\ebar_i} $!)  the
previous formulas give
  $$  \eqalign{
   \exp \left( {\,\log \left( 1 - {_1x_i} \cdot {_2y_i} \right)\, \over
\,{_1x_i}\,} \cdot {_1\hbox{e}_i} \right) \big( {_1L_\lambda} \big)  &  =
{\left( 1 - {_1x_i} \cdot {_2y_i} \right)}^{{\langle \alpha_i \vert \lambda
\rangle} \over 2 \ell}  \cr
   \exp \left( {\,\log \left( 1 - {_1x_i} \cdot {_2y_i} \right)\, \over
\,{_1x_i}\,} \cdot {_1\hbox{e}_i} \right) \big( {_1\fbar_j} \big)  &  =
{_1\fbar_j} -  \cr
   - \delta_{i j} \Bigg( {\,{{\left( 1 - {_1x_i} \cdot {_2y_i} \right)}^{-
1 / \ell} - 1}\, \over \,{_1x_i}\,} \cdot \varepsilon^{d_i}  &
L_{\alpha_i} + {\,{{\left( 1 - {_1x_i} \cdot {_2y_i} \right)}^{1 / \ell} - 1}\, \over \,{_1x_i}\,} \cdot \varepsilon^{-d_i} L_{-\alpha_i} \Bigg) \cdot {_1\ebar_i^{\ell - 1}}  \cr }  $$
and both these are rational functions on  $ \, {\Cal S}_\varepsilon^{(2 \ell)} \times {\Cal S}_\varepsilon^{(2
\ell)} \, $.
                                                   \par
   Now we are left with the case  $ \, x = \ebar_j $,  $ j= \unon \, $.  Consider  $ \, {_1\hbox{e}_i}
\big( {_1\ebar_j} \big) \, $;  if  $ \, \aij = 2 \, $  (i.~e.~$ i=j $)  or
$ \, \aij = 0 \, $  we have  $ \, {_1\hbox{e}_i} \big( {_1\ebar_j} \big) =
0 \, $,  hence
  $$  \exp \left( {\,\log \left( 1 - {_1x_i} \cdot {_2y_i} \right)\, \over
\,{_1x_i}\,} \cdot  {_1\hbox{e}_i} \right) \big( {_1\ebar_j} \big) =
{_1\ebar_j} \; .  $$
   \indent   Therefore we are reduced to make computations in the connected
rank 2 case.  To this end,  we will follow conventions and notations of
[DP], Appendix, and skip for a while left bottom indices "1" (i.~e.~$ \,
{_1\ebar_i} = \ebar_i \, $,  etc.).
                                              \par
   We develop the  $ A_2 $  case; the procedure is the same in the
remaining cases but the computations are longer (cf.~also the  {\it
Remark}  after the proof).
                                              \par
   In this case we have  $ \, d_1 = 1 = d_2 \, $.  Define the root vector
$ \, E_{1 2} := E_{\alpha_1 + \alpha_2} \in \caluqPg \, $  as
  $$  E_{\alpha_{1 2}} \equiv E_{1 2} := T_1(E_2) = - E_1 E_2 + q^{-1} E_2
E_1 \;   \eqno (4.11)  $$ then we have
  $$  E_2 E_1 = q E_1 E_2 + q E_{1 2} \, ,  \qquad \qquad  E_{1 2} E_1 =
q^{-1} E_1 E_{1 2}
\eqno (4.12)  $$
   \indent   Let  $ \Cq(E_1) $  be the field of rational functions in the
indeterminate  $ E_1 $  with  coefficients in  $ \Cq $;  let  $ M $  be
the  $ \Cq(E_1) $--vector  space with  basis  $ \big\{ \ebar_2, \ebar_{1 2}
\big\} $:  then (4.12) tells us that the operation  $ \rho_{E_1} $  of
right multiplication by  $ E_1 $  yields an endomorphism of  $ M $  defined
by the matrix (with respect to the ordered  $
\Cq(E_1) $--basis  $ \, \big\{ \ebar_2, \ebar_{1 2} \big\} \, $)
  $$  \pmatrix
   q E_1  &  0  \\
          q  &  \qm E_1  \\
       \endpmatrix  $$ therefore multiplication by  $ E_1^n $  yields the
endomorphism of  $ M $  defined by the matrix
  $$  {\pmatrix
   q E_1  &  0  \\
       q  &  \qm E_1  \\
      \endpmatrix }^n =
      \pmatrix
   {(q E_1)}^n  &  0  \\
   q {[n]}_{\!q} \cdot E_1^{n-1}  &  {\left( \qm E_1 \right)}^n  \\
      \endpmatrix  \; .   \eqno (4.13)  $$
   \indent   Thus for  $ \, \hbox{e}_1 \left( \ebar_2 \right) \, $  we have
  $$  \displaylines{
   \hbox{e}_1 \left( \ebar_2 \right) := \left[ E_1^{(\ell)},
\ebar_2 \right]{\bigg\vert}_{q=\varepsilon} \!\! =  { E_1^\ell \, \ebar_2 -
\ebar_2 \, E_1^\ell \over  {[\ell]}_{\!q}! }{\Bigg\vert}_{q=\varepsilon}
\!\! =  { E_1^\ell \, \ebar_2 - q^\ell E_1^\ell \, \ebar_2 - q
{[\ell]}_{\!q} \, E_1^{\ell - 1} \, \ebar_{1 2} \over  {[\ell]}_{\!q}!
}{\Bigg\vert}_{q=\varepsilon}
\!\! = \cr
   = { \left( {1 - q^\ell \over q^\ell - q^{-\ell}} \cdot {q-\qm \over
{[\ell - 1]}_{\!q}!} \cdot  E_1^\ell \, \ebar_2 \right)
}{\Bigg\vert}_{q=\varepsilon} - {\left( {q \over {[\ell - 1]}_{\!q}!} \cdot
E_1^{\ell - 1} \, \ebar_{1 2} \right)}{\Bigg\vert}_{q=\varepsilon} =  \cr
   = - {\,\ebar_1^\ell\, \over \,2 \ell\,} \cdot \ebar_2 - {\,\varepsilon\,
\over \,\ell\,} \,
\ebar_1^{\ell - 1} \cdot  \ebar_{1 2} = - {\,x_1\, \over \,2 \ell\,} \cdot
\ebar_2 - {\,\varepsilon\,
\over \,\ell\,} \, \ebar_1^{\ell - 1} \cdot  \ebar_{1 2}  \cr }  $$
(because  $ \, {1 \over {[\ell - 1]}_{\! q}!}{\Big\vert}_{q=\varepsilon} =
{{(q-\qm)}^{\ell - 1} \over
\prod_{s=0}^{\ell - 1} (q^s - q^{-s})}{\Big\vert}_{q=\varepsilon} =
{{(\varepsilon -
\varepsilon^{-1})}^{\ell - 1} \over \ell} \, $);  on the other hand, for  $
\, \hbox{e}_1 \left(
\ebar_{1 2} \right) \, $,  (4.12) gives
  $$  \displaylines{
   \hbox{e}_1(\ebar_{1 2}) := \left[ E_1^{(\ell)}, \ebar_{1 2}
\right]{\bigg\vert}_{q=\varepsilon} =  { E_1^\ell \, \ebar_{1 2} - \ebar_{1
2} \, E_1^\ell \over {[\ell]}_{\!q}! }{\Bigg\vert}_{q=\varepsilon} = {
E_1^\ell \, \ebar_2 - q^{-\ell} E_1^\ell \, \ebar_{1 2} \over
{[\ell]}_{\!q}! }{\Bigg\vert}_{q=\varepsilon} =  \cr
   = {\left( {1 - q^{-\ell} \over q^\ell - q^{-\ell}} \cdot {q-\qm \over
{[\ell - 1]}_{\!q}!} \cdot  E_1^\ell \, \ebar_{1 2}
\right)}{\Bigg\vert}_{q=\varepsilon} = {\,\ebar_1^\ell\, \over \,2 \ell\,}
\cdot
\ebar_{1 2} = {\,x_1\, \over \,2 \ell\,} \cdot \ebar_{1 2}  \cr }  $$
therefore we conclude that  $ \hbox{e}_1 $  restricts to an endomorphism
of  $ M $  defined by the  matrix
  $$  \pmatrix
   -{\ebar_1^\ell \over 2 \ell}  &  0  \\
   -{\varepsilon \over \ell} \cdot \ebar_1^{\ell - 1}  &  {\ebar_1^\ell
\over 2 \ell}  \\
      \endpmatrix =
      \pmatrix
   -{x_1 \over 2 \ell}  &  0  \\
   -{\varepsilon \over \ell} \cdot \ebar_1^{\ell - 1}  &  {x_1 \over 2
\ell}  \\
      \endpmatrix  $$ hence  $ \, \hbox{e}_1^n{\Big\vert}_M = {\left(
\hbox{e}_1{\Big\vert}_M \right)}^n \, $  is given by the matrix
  $$  {\pmatrix
   -{x_1 \over 2 \ell}  &  0  \\
   -{\varepsilon \over \ell} \cdot \ebar_1^{\ell - 1}  &  {x_1 \over 2
\ell}  \\
      \endpmatrix }^n =
      \pmatrix
   {\left( -{x_1 \over 2 \ell} \right)}^n  &  0  \\
   - \delta_{n \in (2 \N + 1)} \cdot {2 \varepsilon \over x_1} \cdot
{\left( {x_1 \over 2 \ell}
\right)}^n \cdot \ebar_1^{\ell - 1}  &  {\left( {x_1 \over 2 \ell}
\right)}^n  \\
      \endpmatrix  $$ (for all  $ \, n \in \N \, $,  where  $ \, \delta_{x
\in Y} := 1 \, $  for  $ \, x \in Y \, $  and
$ \, \delta_{x \in Y} := 0 \, $  for  $ \, x \notin Y \, $),  so that
  $$  {\exp (t \cdot \hbox{e}_1)}{\Big\vert}_M =
   \pmatrix
      e^{-{t \over 2 \ell} \cdot x_1}  &  0  \\
      -{\varepsilon \over x_1} \cdot \left(e^{{t \over 2 \ell} \cdot x_1} -
e^{-{t \over 2 \ell} \cdot  x_1} \right) \cdot \ebar_1^{\ell - 1}  &  e^{{t
\over 2 \ell} \cdot x_1}  \\
   \endpmatrix  $$ where  $ t $  denotes any indeterminate which commute
with
$ E_1 $;  in particular for  $ \, t = {\log(1 - w_1) \over \ebar_1^\ell} \,
$,  with  $ \, w_1 := {_1\ebar_1}^\ell \cdot {_2\fbar_1}^\ell = {_1x_1}
\cdot {_2y_1} \, $,  we get
  $$  {\exp \left( {\log(1 - w_1) \over \ebar_1^\ell} \cdot \hbox{e}_1
\right)}{\Bigg\vert}_M =
   \pmatrix
      {(1 - w_1)}^{-{1 \over 2 \ell}}  &  0  \\
      -{\varepsilon \over \ebar_1^\ell} \cdot \left({(1 - w_1)}^{1 \over 2
\ell} - {(1 - w_1)}^{-{1
\over 2 \ell}} \right) \cdot \ebar_1^{\ell - 1}  &  {(1 - w_1)}^{1 \over 2
\ell}  \\
   \endpmatrix  $$
thus
  $$  \displaylines{
   \exp \left( {\log(1 - {_1x_1} \cdot {_2y_1}) \over {_1x_1}} \cdot
{_1\hbox{e}_1} \right) \left(  {_1\ebar_2} \right) =  \cr
   = {\left( 1 - {_1x_1} \cdot {_2y_1} \right)}^{-{1 \over 2 \ell}} \cdot
{_1\ebar_2} - {\varepsilon \over {_1x_1}} \cdot \left({ \left( 1 - {_1x_1}
\cdot {_2y_1} \right)}^{1 \over 2 \ell} -  {\left( 1 - {_1x_1} \cdot
{_2y_1} \right)}^{-{1 \over 2 \ell}} \right) \cdot {_1\ebar_1}^{\ell - 1}
\cdot {_1\ebar_{1 2}}  \cr }  $$  which is a rational section of a  $ \,
{\Cal S}_\varepsilon^{(2 \ell)} \times  {\Cal S}_\varepsilon^{(2 \ell)} \,
$,  q.~e.~d.
                                                      \par
   As for  $ \, \exp \left( {\log \left( 1 - {_1x_2} \cdot {_2y_2} \right)
\over {_1x_2}} \cdot  {_1\hbox{e}_2} \right) \, $,  everything comes from
above by symmetry, namely because  $ \, \alpha_2 = s_1 s_2 (\alpha_1) \, $
implies  $ \, {_1\ebar_2} = T_1 T_2 \left( {_1\ebar_1} \right) \, $,  $ \,
{_1\fbar_2} = T_1 T_2 \left( {_1\fbar_1} \right) \, $,  and  $ \,
{_1\hbox{e}_2} = (T_1 T_2) \smallcirc {_1\hbox{e}_1} \smallcirc {(T_1
T_2)}^{-1} \, $;  on the other hand, in the other cases of rank 2 (that is
$ B_2 $  and  $ G_2 \, $)  such a symmetric situation does not occur, hence
we must perform direct computation for  $ \, \exp \! \left( {\log \left( 1
- {_1x_2} \cdot {_2y_2} \right) \over {_1x_2}}
\cdot {_1\hbox{e}_2} \right) \, $  too (this is entirely similar, although
longer, to the previous one).
                                                \par
   Finally, it is clear that restricting to subalgebras
$ Z_\varepsilon $  and  $ Z_0 $  we get (bi)rational Poisson automorphisms
of their spectra, by the  same argument of the end of the proof of
Proposition 4.2.   $ \square $
\enddemo

\vskip7pt

{\it Remark:} \; the very (theoretical) reason why computations do work in
{\it all}  rank two cases, so that Theorem 4.3 does hold, lies in the
availability of the commutation formulas for quantum root vectors (the
so-called Levendorskij-Soibel'man formulas, cf.~[DP], Theorem 9.3),
strictly related with the existence of a convex ordering on the set of
positive roots.

\vskip7pt

   The previous result can be still improved when considering the central
Hopf subalgebra  $ Z_0 $,   hence the Poisson group  $ H $,  as the
following shows:

\vskip7pt

\proclaim{Theorem 4.4}  The birational Poisson automorphism  $ \, {\Cal
R}_{\varepsilon, \ell}^{\,\ast}
\colon {\Cal H}^{(2 \ell)} \times {\Cal H}^{(2 \ell)} \rightarrow {\Cal
H}^{(2 \ell)} \times {\Cal H}^{(2 \ell)}
\, $  pushes down to a birational Poisson automorphism  $ \, {\Cal
R}_{\varepsilon, \ell}^{\,\ast} \colon {\Cal H}^{(2)} \times {\Cal H}^{(2)}
\rightarrow {\Cal H}^{(2)} \times {\Cal H}^{(2)} \, $,  independent of
$ \ell $,  of a  $ 2 $-fold  ramified covering  $ \, {\Cal H}^{(2)} \times
{\Cal H}^{(2)} \, $  of  $ \, H \times
H \, $;  moreover,  $ \, {\Cal R}_{\varepsilon,
\ell}^{\,\ast} \neq \tau^{(2)} $  (the "twist" map of  $ \, {\Cal H}^{(2)}
\times {\Cal H}^{(2)}
\, $),  and  $ {\Cal R}_{\varepsilon, \ell}^{\,\ast} $  enjoys the dual
properties of (1.4--6):  in particular,
$ \, m \left( {\Cal R}_{\varepsilon, \ell}^{\,\ast} (x,y)
\right) = y \cdot x \, $  for all  $ \, x, y \in {\Cal H}^{(2)} \, $
(where  $ m $  and  "$ \, \cdot \, $"  denote the product of
$ \, {\Cal H}^{(2)} \, $),  and a braid group action exists on
$ \times $--powers  of  $ {\Cal H}^{(2)} $.
\endproclaim

\demo{Proof}  As for Theorem 4.3, the proof amounts to check that some
series do converge on an appropriate covering.  Namely, we have to check
that
  $$  \exp \left( {\, - \varepsilon \, \over \, 2 d_\alpha
\ell^{\scriptscriptstyle 2} \,} \, \psi \!\left( x_\alpha
\otimes y_\alpha \right) \cdot \hbox{ad}_{\{\ ,\ \}} \!\left( x_\alpha
\otimes y_\alpha \right) \right) ({_1w} \otimes {_2w})  $$ does converge to
a rational function on a covering  $ \, {\Cal H}^{(2)} \times {\Cal
H}^{(2)} \, $  as claimed for all  $ \, \alpha \in R^+ \, $  and for all  $
\, {_iw} \in \big\{\, 1, {_ix_\beta}, {_iz_\lambda}, {_iy_\gamma}
\,\big\vert\, \beta, \gamma \in R^+; \, \lambda \in P \,\big\} \, $,  $ \,
i = 1, 2 \, $.  This again amounts to prove the same for functions
  $$  \eqalign{
   \exp  &  \left( {\log ( 1 - {_1x_\alpha} \cdot {_2y_\alpha} ) \over
{_1x_\alpha}} \cdot {_1\hbox{e}_\alpha} \right) ({_1w})  \cr
   \exp  &  \left( {\log ( 1 - {_1x_\alpha} \cdot {_2y_\alpha} ) \over
{_2y_\alpha}} \cdot {_2\hbox{f}_\alpha} \right) ({_2w})  \cr }  $$ for all
$ \alpha $  and  $ {_iw} $  like above.  As for Theorem 4.3, we deal with
the first function, the proof for the second one following by symmetry.
                                                   \par
   By the braid group action we can again reduce to the case of simple
roots  $ \, \alpha =
\alpha_i \, $.  Furthermore (cf.~[DK], \S 3.4, and [DP], \S 19), with
respect to coordinates  $ \, x_\gamma := \ebar_\gamma^\ell \, $,  $ \,
z_\lambda := L_\lambda^\ell \, $,  $ \, y_\gamma :=
\fbar_\gamma^\ell \, $,  {\sl the formulas for derivations}  $
\hbox{e}_\alpha $  {\sl are independent of}  $ \ell \, $:  therefore we can
fix  $ \, \ell = 1 \, $  and perform computations in  $ U_1 $.
                                                   \par
   Again direct computation (or formulas in the proof of [DK], Proposition
3.5) gives
  $$  \exp \left( t \cdot {_1\hbox{e}_i} \right) ({_1z_\lambda}) = e^{
(\langle \alpha_i \vert \lambda \rangle / 2) \cdot t \cdot {_1\!x_i}} \cdot
{_1z_\lambda}   \eqno (4.14)  $$ for any indeterminate  $ t $  which
commutes with
$ \, {_1\ebar_i} = {_1x_i} \, $;  then for  $ \, t = {\log ( 1 - {_1x_i}
\cdot {_2y_i} ) \over {_1x_i}} \, $  we have
  $$  \exp \left( {\log ( 1 - {_1x_i} \cdot {_2y_i} ) \over {_1x_i}} \cdot
{_1\hbox{e}_i} \right) ({_1z_\lambda}) = {\left( 1 - {_1x_i} \cdot {_2y_i}
\right)}^{{\langle \alpha_i \vert \lambda \rangle} \over 2} \cdot
{_1z_\lambda}  $$  which is a rational function on a 2-fold ramified
covering  $ \, {\Cal H}^{(2)} \times {\Cal H}^{(2)} \, $  of
$ \, H \times H \, $.
                                                      \par
   Now consider  $ \, \exp \left( {\log ( 1 - {_1x_i} \cdot {_2y_i} ) \over
{_1x_i}} \cdot  {_1\hbox{e}_i} \right) \left( {_1x_\gamma} \right) \, $,
with  $ \, \gamma \in R^+ \, $  (notice that now simple root vectors  $ \,
\ebar_j = x_j \, $  ($ j= \unon $)  are not enough to generate  $ U_1^+ $
(the "positive part" of  $ U_1 $):  we do need all root vectors  $ \,
\ebar_\gamma = x_\gamma \, $,  $ \,
\gamma \in R^+ \, $).  For any fixed pair  $ (\alpha, \gamma) $  of
positive roots, let us denote by  $ R^+_{\beta, \gamma} $  the rank 2 root
system spanned by  $ \{\alpha, \gamma\} $  in  $ R^+ $.  The following is
well known (cf.~e.~g.~[DP], first Lemma of \S 15.4):
                                                       \par
   {\it Claim:} \  For any fixed pair  $ (\alpha_i,
\gamma) $  of positive roots with  $ \alpha_i $  simple, there exists  $ w
\in W $  and  $ \alpha_1, \alpha_2 \in R^+ $  such that  $ \, w \left(
R^+_{\alpha_1, \alpha_2} \right) = R^+_{\alpha_i, \gamma} \, $  and
$ \, w(\alpha_1) = \alpha_i \, $.
                                                      \par
   Thanks to  {\it Claim\/}  and (4.10) we are reduced to make computations
in the rank 2 case; the same  holds when considering  $ \, \exp \left(
{\log ( 1 - {_1x_i} \cdot {_2y_i} ) \over {_1x_i}} \cdot {_1\hbox{e}_i}
\right) \left( {_1y_\gamma} \right) \, $,  with  $ \, \gamma \in R^+ \, $
(now again negative simple root vectors  $ \, \fbar_j = y_j \, $  ($ j=
\unon $)  are not enough to generate  $ U_1^- $  (the "negative part" of  $
U_1 $):  we do need all negative root vectors  $ \, \fbar_\gamma = y_\gamma
\, $,  $ \, \gamma \in R^+ \, $).  We denote by  $ T $  the type of a root
system of rank 2 (hence  $ \, T \in \big\{ A_1 \times A_1, A_2, B_2, G_2
\big\} \, $).
                                                      \par
   $ T = A_1 \times A_1 $: \  First of all, since  $ \, \hbox{e}_j \left(
x_j \right) = \hbox{e}_j \left( \ebar_j \right) = 0 \, $  ($ j= 1, 2 $),
we have
  $$  \exp \left( {\log ( 1 - {_1x_j} \cdot {_2y_j} ) \over {_1x_j}} \cdot
{_1\hbox{e}_j} \right) \left( {_1x_j} \right) = {_1x_j}  \qquad \qquad  (j=
1, 2)  \; ;  $$ second, since  $ a_{1 2} = 0 $,  we have  $ \, \hbox{e}_i
\left( x_j \right) = 0 \, $  (for  $ \, i, j \in \{1, 2\}
\, $,  $ i \not= j $)  whence
  $$  \exp \left( {\log ( 1 - {_1x_i} \cdot {_2y_i} ) \over {_1x_i}} \cdot
{_1\hbox{e}_i} \right) \left( {_1x_j} \right) = {_1x_j}  $$ (for  $ \, i, j
\in \{1, 2\} \, $,  $ i \not= j $) thus we are done with generators  $
x_\alpha $'s.
                                                      \par
   As for negative root vectors  $ \, y_\alpha = \fbar_\alpha \, $,  we
have  $ \, \hbox{e}_i (y_j) =
\delta_{i j}  \cdot (z_{\alpha_i} - z_{-\alpha_i}) \, $,  whence
  $$  {\hbox{e}_i}^{\!n} (y_j) = \delta_{i j} \cdot \left(
{\hbox{e}_i}^{\!n-1} (z_{\alpha_i}) - {\hbox{e}_i}^{\!n-1} (z_{-\alpha_i})
\right) = \delta_{i j} \cdot \left( {(-x_i)}^{n-1} \cdot z_{\alpha_i} -
x_i^{\,n-1} \cdot z_{-\alpha_i} \right)  $$
(thanks to (4.14)) for all  $ \, n \in \N_+ \, $,  thus
  $$  \exp \left( t \cdot {_1\hbox{e}_i} \right) ({_1y_j}) = {_1y_j} -
\delta_{i j} \cdot \left( {e^{- t \cdot {_1x_i}} - 1 \over {_1x_i}} \cdot
z_{\alpha_i} + {e^{t \cdot {_1x_i}} - 1 \over {_1x_i}} \cdot z_{-\alpha_i}
\right)  $$ for any indeterminate  $ t $  which commutes with
$ {_1x_i} $;  for  $ \, t = {\log ( 1 - {_1x_i} \cdot {_2y_i} ) \over
{_1x_i}} \, $  we get
  $$  \displaylines{
   \exp \left( {\log ( 1 - {_1x_i} \cdot {_2y_i} ) \over {_1x_i}} \cdot
{_1\hbox{e}_i} \right) ({_1y_j}) =  \cr
   = {_1y_j} - \delta_{i j} \cdot \left( { {\left( 1 - {_1x_i} \cdot
{_2y_i} \right)}^{-1} - 1 \over {_1x_i} } \cdot z_{\alpha_i} + { \left( 1 -
{_1x_i} \cdot {_2y_i} \right) - 1 \over {_1x_i} } \cdot z_{-\alpha_i}
\right)  \cr }  $$
which is a rational function on a 2-fold ramified covering  $ \,
{\Cal H}^{(2)} \times {\Cal H}^{(2)} \, $  of  $ \, H \times H \, (\, =
Spec(Z_0) \times Spec(Z_0)
\,) \, $.  Since for  $ \, T = A_1 \times A_1 \, $  we have  $ \, R^+ =
\{\alpha_1, \alpha_2\} \, $, we are done.
                                                      \par
   $ T = A_2 $: \  We follow again conventions and notations of [DP],
Appendix.  In the present case we  have  $ \, d_1 = 1 = d_2 \, $,  and  $
\, R^+ = \{\alpha_1, \alpha_{1 2}:= \alpha_1 + \alpha_2,
\alpha_2\} \, $,  and we define the root vector  $ \, E_{1 2} := - E_1 E_2 + \qm E_2 E_1 \, $  (cf.~(4.11)).  For  $ \, \gamma = \alpha_1 \, $  we have as above
  $$  \exp \left( {\log ( 1 - {_1x_1} \cdot {_2y_1} ) \over {_1x_1}} \cdot {_1\hbox{e}_1} \right)  ({_1x_1}) = {_1x_1}  \; .  $$
   \indent   Then let  $ M $  be the  $ \Cq(E_1) $--vector space  with basis  $ \, \left\{ x_2, x_{1 2} \right\} = \left\{ \ebar_2, \ebar_{1 2} \right\} \, $:  then (4.12) tells us that the operation of right multiplication by
$ E_1 $  yields an endomorphism of  $ M $  defined by the matrix (with respect to  $ \, \left\{ \ebar_2, \ebar_{1 2} \right\} \, $)
  $$  \pmatrix
   q E_1  &  0  \\
       q  &  \qm E_1  \\
      \endpmatrix  \, .  $$
   \indent   Thus for  $ \, \hbox{e}_1(x_2) \, $  we  have
  $$  \displaylines{
   \hbox{e}_1(x_2) := \left[ E_1, \ebar_2 \right]{\Big\vert}_{q=1} = \left( { E_1
\ebar_2 - \ebar_2 E_1 } \right) {\Big\vert}_{q=1} = \left( { E_1 \ebar_2 - q E_1 \ebar_2 -
q \ebar_{1 2} } \right) {\Big\vert}_{q=1} =  \cr
   = {\left( {1 - q \over q - \qm} \cdot \ebar_1 \ebar_2
\right)}{\Bigg\vert}_{q=1} - {\left( q \cdot
\ebar_{1 2} \right)}{\Big\vert}_{q=1} = - {\,\ebar_1\, \over \,2\,} \cdot
\ebar_2 - \ebar_{1 2} =  - {\,x_1\, \over \,2\,} \cdot x_2 - x_{1 2} \; ;
\cr }  $$  on the other hand, for  $ \, \hbox{e}_1(x_{1 2}) \, $,  (4.12)
gives
  $$  \displaylines{
   \hbox{e}_1(x_{1 2}) := \left[ E_1, \ebar_{1 2} \right]{\Big\vert}_{q=1}
= \left( { E_1 \ebar_{1 2} -
\ebar_{1 2} E_1 } \right) {\Big\vert}_{q=1} = \left( { E_1 \ebar_2 - \qm E_1 \ebar_{1 2} } \right) {\Big\vert}_{q=1} =  \cr
   = {\left( {1 - \qm \over q - \qm} \cdot E_1 \ebar_{1 2}
\right)}{\Bigg\vert}_{q=1} = {\,\ebar_1\,
\over \,2\,} \cdot \ebar_{1 2} = {\,x_1\, \over \,2\,} \cdot x_{1 2}  \cr
}  $$  therefore we conclude that  $ \hbox{e}_1 $  restricts to an
endomorphism of  $ M $  defined by the  matrix
  $$  \pmatrix
   -{x_1 \over 2}  &  0  \\
   -1  &  {x_1 \over 2}  \\
      \endpmatrix  $$ hence  $ \, \hbox{e}_1^n{\Big\vert}_M = {\left(
\hbox{e}_1{\Big\vert}_M \right)}^n \, $  is given by the matrix
  $$  {\pmatrix
   -{x_1 \over 2}  &  0  \\
   -1  &  {x_1 \over 2}  \\
      \endpmatrix }^n =
     \pmatrix
   {\left( -{x_1 \over 2} \right)}^n  &  0  \\
   - \delta_{n \in (2 \N + 1)} \cdot {2 \over x_1} \cdot {\left( {x_1 \over
2} \right)}^n  &  {\left( {x_1 \over 2} \right)}^n  \\
      \endpmatrix  $$ (for all  $ \, n \in \N \, $,  so that
  $$  {\exp (t \cdot \hbox{e}_1)}{\Big\vert}_M =
   \pmatrix
      e^{-{t \over 2} \cdot x_1}  &  0  \\
      -{1 \over x_1} \cdot \left(e^{{t \over 2} \cdot x_1} - e^{-{t \over
2} \cdot x_1} \right)  &  e^{{t \over 2} \cdot x_1}  \\
   \endpmatrix  $$
for any  $ t $  which commutes with
$ E_1 $,  and for  $ \, t = {\log(1 - {_1x_1} \cdot {_2y_1}) \over {_1x_1}} \, $  we have
  $$  \displaylines{
   {\exp \left( {\log(1 - {_1x_1} \cdot {_2y_1}) \over x_1} \cdot {_1\hbox{e}_1} \right)}{\Bigg\vert}_M =  \cr
   = \pmatrix
      {(1 - {_1x_1} \cdot {_2y_1})}^{-{1 \over 2}}  &  0  \\
      -{1 \over {_1x_1}} \cdot \left({(1 - {_1x_1} \cdot {_2y_1})}^{1 \over
2} - {(1 - {_1x_1} \cdot {_2y_1})}^{-{1 \over 2}} \right)  &  {(1 - {_1x_1}
\cdot {_2y_1})}^{1 \over 2}  \\
   \endpmatrix   \cr }  $$ or, in other words,
  $$  \displaylines{
   \exp \left( {\log(1 - {_1x_1} \cdot {_2y_1}) \over {_1x_1}} \cdot
{_1\hbox{e}_1} \right) ({_1x_2}) =  \cr
   = {(1 - {_1x_1} \cdot {_2y_1})}^{-{1 \over 2}} \cdot {_1x_2}  - {1 \over
{_1x_1}} \cdot \left({(1 - {_1x_1} \cdot {_2y_1})}^{1 \over 2} - {(1 -
{_1x_1} \cdot {_2y_1})}^{-{1 \over 2}} \right) \cdot {_1x_{1 2}}  \cr
   \exp \left( {\log(1 - {_1x_1} \cdot {_2y_1}) \over {_1x_1}} \cdot
{_1\hbox{e}_1} \right) ({_1x_{1 2}}) = {(1 - {_1x_1} \cdot {_2y_1})}^{1
\over 2} \cdot {_1x_{1 2}}
\cr }  $$  and these are rational functions on a  $ 2 $-fold  ramified
covering  $ {\Cal H}^{(2)} \times {\Cal H}^{(2)} $  of  $ H \times H $,
q.~e.~d.
                                                      \par
   For negative root vectors  $ \, y_\alpha = \fbar_\alpha $'s,  define  $
\, F_{\alpha_{1 2}} \equiv F_{1 2} := T_1(F_2) = - F_2 F_1 + q F_1 F_2 \,
$;  then we have again  $ \, \hbox{e}_1 (y_j) = \delta_{1 j} \cdot
(z_{\alpha_1} - z_{-\alpha_1}) \, $  ($ j= 1, 2 $),  whence
  $$  \displaylines{
   \exp \left( {\log ( 1 - {_1x_1} \cdot {_2y_1} ) \over {_1x_1}} \cdot
{_1\hbox{e}_1} \right) ({_1y_j}) =  \cr
   = {_1y_j} - \delta_{1 j} \cdot \left( { {\left( 1 - {_1x_1} \cdot
{_2y_1} \right)}^{-1} - 1 \over {_1x_1} } \cdot z_{\alpha_1} + {\left( 1 -
{_1x_1} \cdot {_2y_1} \right) - 1 \over {_1x_1}} \cdot z_{-\alpha_1}
\right)
\cr }  $$  ($ j= 1, 2 $)  which is a rational function on the proper
covering; this takes care of  $ \gamma =
\alpha_1 $  and  $ \gamma = \alpha_2 $.  At last, for  $ \, \gamma =
\alpha_{1 2} := \alpha_1 + \alpha_2
\, $,  we have
  $$  \displaylines{
   \hbox{e}_1 (y_{1 2}) = {\left[ E_1, \fbar_{1 2}
\right]}{\Big\vert}_{q=1} = {\left( - F_2
\left[ E_1, \fbar_1 \right] + q \left[ E_1, \fbar_1 \right] F_2
\right)}{\Big\vert}_{q=1} =  \cr
   = {\left( - \fbar_2 \left( L_{\alpha_1} - L_{-\alpha_1} \right) + q
\left( L_{\alpha_1} -  L_{-\alpha_1} \right) \fbar_2
\right)}{\Big\vert}_{q=1} = {\left( {q^2 - 1 \over q - \qm} \cdot
\fbar_2 \, L_{\alpha_1} \right)}{\Big\vert}_{q=1} = z_{\alpha_1} \cdot y_2
\cr }  $$ now, since for all  $ \, n \in \N \, $  we have
  $$  \hbox{e}_1^{\,n} (z_{\alpha_1} \cdot y_2) = \hbox{e}_1^{\,n}
(z_{\alpha_1}) \cdot y_2 = {(-x_1)}^n \cdot z_{\alpha_1} \cdot y_2  $$ we
get, for all  $ \, n \in \N_+ \, $
  $$  \hbox{e}_1^{\,n} (y_{1 2}) = \hbox{e}_1^{\,n-1} (z_{\alpha_1} \cdot
y_2) = \hbox{e}_1^{\,n-1} \cdot z_{\alpha_1} y_2 = {(-x_1)}^{n-1} \cdot
z_{\alpha_1} y_2 = -{{(-x_1)}^n \over x_1} \cdot z_{\alpha_1} y_2  $$
whence  $ \, \exp \left( t \cdot \hbox{e}_1 \right) (y_{1 2}) = y_{1 2} -
{e^{-t \cdot x_1} - 1 \over x_1} \cdot z_{\alpha_1} y_2 \, $  and finally
  $$  \exp \left( {\log ( 1 - {_1x_1} \cdot {_2y_1} ) \over {_1x_1}} \cdot
{_1\hbox{e}_1} \right)  ({_1y_{1 2}}) = {_1y_{1 2}} - {{(1 - {_1x_1} \cdot
{_2y_1})}^{-1} - 1 \over {_1x_1}} \cdot {_1z_{\alpha_1}}  \cdot {_1y_2}  $$
the latter being a rational function on the covering  $ \, {\Cal H}^{(2)}
\times {\Cal H}^{(2)} \, $  of
$ \, H \times H \, $,  q.~e.~d.
                                                     \par
   As for  $ \, \exp \left( {\log(1 - {_1x_2} \cdot {_2y_2}) \over {_1x_2}}
\cdot {_1\hbox{e}_2} \right)
\, $,  everything follows by symmetry; on the other hand, in cases  $ B_2
$  and  $ G_2 $  such a symmetric situation does not occur, hence we must
perform direct computation for  $ \, \exp \left( {\log(1 - {_1x_2} \cdot
{_2y_2}) \over {_1x_2}} \cdot {_1\hbox{e}_2} \right) \, $  too (which is
completely similar, although quite longer, to the previous one).   $ \square $
\enddemo

\vskip7pt

   We stress the fact that the proof of Theorem 4.4 above also contains the proof
of the following one, which  means that the adjoint action of the  $ R $-matrix
does specialize for  $ \, q \rightarrow 1 \, $  to something more than formal,
with a very precise geometric meaning:

\vskip7pt

\proclaim{Theorem 4.5}  For the braided Hopf algebra  $ \,
\big( F[H], {\Cal R}_1 \big) = (U_1, {\Cal R}_1) \, $  the "formal
automorphism"  $ {\Cal R}_1 $  is in fact an effective Poisson automorphism
of the field of rational functions on  $ {\Cal H}^{(2)} \times {\Cal
H}^{(2)} $.  In other words,  $ {\Cal R}_1 $  defines a Poisson birational
automorphism  $ {\Cal R}_1^{\,\ast} $  of
$ {\Cal H}^{(2)} \times {\Cal H}^{(2)} $  which enjoys the dual properties
of (1.4--6), in particular  $ \, m \left( {\Cal R}_1^{\,\ast} (x,y) \right)
= y \cdot x \, $  for all  $ \, x, y \in {\Cal H}^{(2)} \, $,  and a braid
group action exists on  $ \times $--powers  of  $ {\Cal H}^{(2)} $.   $ \square $
\endproclaim

\vskip7pt

{\bf 4.6.} \;  Recall that, by general theory, one has  $ \, Spec(A \otimes B) =
Spec(A) \times Spec(B) \, $  for all associative unital algebras  $ A $  and  $ B $;  moreover,
if  $ M $  and  $ N $  are Poisson manifolds then  $ \, M \times N \, $  is a Poisson manifold too, whose symplectic leaves are all the products of symplectic leaves of  $ M $  and  $ N $.  In our context this implies that the symplectic leaves of  $ \, H^{\times 2} \, $  (resp.~$ {\left( {\Cal H}^{(N)} \right)}^{\times 2} \, $)  are all the products of symplectic leaves of  $ H $  (resp.~$ {\Cal H}^{(N)} \, $).
                                                    \par
   Let  $ \, N \in \{2 \ell, \infty \} $.  Let us denote by  $ {\Cal Z}_0 $  the algebra of meromorphic functions on
$ {\Cal H}^{(N)} $.  As we said, we can look at
$ S_\varepsilon $  as a sheaf of algebras over  $ H \, $;  similarly, we can look at  $ {\Cal S}_\varepsilon^{(N)} $  (the  $ N $--fold  covering of  $ S_\varepsilon \, $)  as a sheaf of algebras over  $ {\Cal H}^{(N)} \, $:  its algebra of meromorphic sections  $ {\Cal U}_\varepsilon $  is then  $ \, {\Cal U}_\varepsilon = {\Cal Z}_0 \otimes_{Z_0} U_\varepsilon $.  Now represent the elements of  $ \, {\Cal H}^{(N)} = Spec({\Cal Z}_0) \, $  as maximal ideals of  $ {\Cal Z}_0 $,  and let  $ \, {\frak m} \in {\Cal H}^{(N)} \, $:  the fibre over  $ \frak m $  of our sheaf is then  $ \, {\Cal U}_\varepsilon \big/ {\frak m} \, {\Cal U}_\varepsilon \, $.  Similarly, the fibre over  $ \, ({\frak m}, {\frak n}) \in {\Cal H}^{(N)} \times {\Cal H}^{(N)} \, $  (of the sheaf of algebras  $ {\Cal S}_\varepsilon^{(N)} \times {\Cal S}_\varepsilon^{(N)} $  over  $ {\Cal H}^{(N)} \times {\Cal H}^{(N)} $)  is  $ \, ({\Cal U}_\varepsilon \otimes {\Cal U}_\varepsilon) \big/ ({\frak m} \, {\Cal U}_\varepsilon \otimes {\Cal U}_\varepsilon + {\Cal U}_\varepsilon \otimes {\frak n} \, {\Cal U}_\varepsilon) \, $.

\vskip7pt

\proclaim {Proposition 4.7}  $ {\Cal R}_{\varepsilon,N}^* $  is a meromorphic automorphism of  $ {\Cal S}_\varepsilon^{(N)} \times {\Cal S}_\varepsilon^{(N)} $  as a fibre bundle over  $ {\Cal H}^{(N)} \times {\Cal
H}^{(N)} $  with respect to the meromorphic automorphism
$ {\Cal R}_{\varepsilon,N}^\ast {\Big\vert}_{{\Cal H}^{(N)} \times {\Cal H}^{(N)}} $  of the base space  $ {\Cal H}^{(N)} \times {\Cal H}^{(N)} $;  in other words, the following diagram is commutative
  $$  \CD
       {\Cal S}_\varepsilon^{(N)} \times {\Cal S}_\varepsilon^{(N)}  @>{{\Cal R}_{\varepsilon,N}^\ast}>>   {\Cal S}_\varepsilon^{(N)} \times {\Cal S}_\varepsilon^{(N)}  \\
       @V{\pi}VV                           @VV{\pi}V  \\
       {\Cal H}^{(N)} \times {\Cal H}^{(N)}   @>>{{\Cal
R}_{\varepsilon,N}^\ast}>   {\Cal H}^{(N)} \times {\Cal H}^{(N)}  \\
      \endCD  $$
where  $ \, \pi \colon \, {\Cal S}_\varepsilon^{(N)}
\times {\Cal S}_\varepsilon^{(N)} \longrightarrow {\Cal H}^{(N)} \times {\Cal H}^{(N)} \, $  is the projection map of the fibre bundle.  In particular  $ {\Cal R}_{\varepsilon,N}^* $  leaves invariant the fibres of
$ {\Cal S}_\varepsilon^{(N)} \times {\Cal S}_\varepsilon^{(N)} $  over symplectic leaves of  $ {\Cal H}^{(N)} \times {\Cal H}^{(N)} $  (i.~e.~the
preimages, with respect to  $ \pi $,  of symplectic leaves of  $ {\Cal H}^{(N)} \times {\Cal H}^{(N)} $).
\endproclaim

\demo{Proof}  This is more or less trivial, by construction.  Let  $ \,
\frak s = (\frak m, \frak n) \in Spec \left( {{\Cal Z}_0}^{\otimes 2}
\right) = {\Cal H}^{(N)} \times {\Cal H}^{(N)} \, $  be a maximal ideal of
$ {{\Cal Z}_0}^{\otimes 2} $;  then its fibre  $ {{\Cal
U}_\varepsilon}^{\otimes 2} \Big/ {\frak s} \, {{\Cal
U}_\varepsilon}^{\otimes 2} $  is mapped by  $ {\Cal R}_{\varepsilon,N}^*
$  onto  $ \, {\Cal R}_{\varepsilon,N}^\ast \left( {{\Cal
U}_\varepsilon}^{\otimes 2} \Big/ {\frak s} \, {{\Cal
U}_\varepsilon}^{\otimes 2} \right) = {{\Cal U}_\varepsilon}^{\otimes 2}
\Big/ {{\Cal R}_\varepsilon}^{\! -1}({\frak s}) \, {{\Cal
U}_\varepsilon}^{\otimes 2} \, $,  whence everything easily follows.   $ \square $
\enddemo

\vskip3truecm

\Refs
\endRefs

\vskip8pt

\smallrm

[DD] \  I.~Damiani, C.~De Concini, {\smallit Quantum groups and Poisson
groups,\/}  in W.~Baldoni,  M.~Picardello (eds.), {\smallit Representations
of Lie groups and quantum groups,\/} Longman Scientific
$ \! {\scriptstyle \and} \! $  Technical.

\vskip5pt

[DK] \  C.~De Concini, V.~G.~Kac, {\smallit Representations of Quantum
Groups at Roots of 1,\/} in  Colloque Dixmier 1989, Progr.~in
Math.~{\smallbf 92} (1990), 471 -- 506.

\vskip5pt

[DKP] \  C.~De Concini, V.~G.~Kac, C.~Procesi, {\smallit  Quantum coadjoint
action,\/} Jour.~Am. Math.~Soc.~{\smallbf 5} (1992), 151 -- 189.

\vskip5pt

[DL] \  C.~De Concini, V.~Lyubashenko, {\smallit Quantum function algebra
at roots of 1,\/}  Adv.~in  Math.~{\smallbf 108} (1994), 205--262.

\vskip5pt

[DP] \  C.~De Concini, C.~Procesi, {\smallit Quantum Groups,\/}  in
L.~Boutet de Monvel, C.~De Concini,  C.~Procesi, P.~Schapira, M.~Vergne
(eds.),  {\smallit D-modules, Representation Theory, and Quantum
Groups,\/}  Lectures Notes in Mathematics  {\smallbf 1565},  Springer  $ \!
{\scriptstyle \and} \! $  Verlag,  Berlin-Heidelberg-New York, 1993.

\vskip5pt

[Dr] \  V.~G.~Drinfeld, {\smallit  Quantum Groups,\/} in {\smallit
Proc.~Intern.~Congress of Math.\/}  (Berkeley, 1986), AMS 1987, pp.
798--820.

\vskip5pt

[Ex] \  H.~Exton, {\smallit  q--Hypergeometric  Functions and
Applications,\/}  Ellis  Hordwood Series  {\smallit Mathematics and its
Applications\/},  Ellis Hordwood Limited 1983.

\vskip5pt

[GR] \  G.~Gasper, M.~Rahman, {\smallit Basic hypergeometric  series,\/}
Encyclopedia of Mathematics and its Applications {\smallbf 35}, Cambridge
University Press,  1990.

\vskip5pt

[Ji] \  M.~Jimbo, {\smallit  A q--difference analogue of  U${\scriptstyle
({\frak g})}$   and  the Yang-Baxter equation,\/}  Lett.~Math.
Phys.~{\smallbf 10} (1985), 63--69.

\vskip5pt

[KR] \  A.~N.~Kirillov, N.~Reshetikin, {\smallit  q--Weyl Group and a
Multiplicative  Formula for Universal R--Matrices,\/}
Comm.~Math.~Phys.~{\smallbf 134} (1990), 421--431.

\vskip5pt

[LS] \  S.~Z.~Levendorskij, Ya.~S.~Soibel'man,  {\smallit Some applications
of the quantum Weyl  groups,\/}  J.~Geom. Physics {\smallbf 7} (1990),
n$^{\hbox{\ssmallrm o}}$ 2, 241--254.

\vskip5pt

[Re] \  N.~Reshetikin, {\smallit  Quasitriangularity of quantum groups at
roots of 1,\/}  Commun.~Math.~Phys.~{\smallbf 170} (1995), 79--99.

\vskip5pt

[Ta] \  T.~Tanisaki, {\smallit  Killing forms, Harish-Chandra
Isomorphisms, and Universal R--Matrices for Quantum Algebras,\/}
Internat.~J.~Modern Phys.~A {\smallbf 7},  Suppl.~1 B (1992), 941--961.

\vskip1truecm

\enddocument

\bye
\end